\documentclass[]{interact}

\usepackage{amsmath}
\usepackage{amssymb}
\usepackage{amsthm}
\usepackage{mathrsfs}
\usepackage{amsfonts}
\usepackage{upgreek}
\usepackage{dsfont}

\usepackage{cite}

\usepackage{array}
\usepackage{float}
\usepackage[caption=false]{subfig}
\usepackage{bm}

\usepackage{url}

\usepackage{booktabs}
\usepackage{multirow}

\usepackage{rotating}

\usepackage{epstopdf}
\usepackage{placeins}

\numberwithin{equation}{section}

\DeclareMathOperator*{\argmin}{argmin}
\DeclareMathOperator*{\diag}{diag}

\newcommand{\de}{\,\mathrm{d}}

\newcommand{\bigo}{\mathcal{O}}

\newcommand*{\lvec}[1]{\mathbf{#1}}
\newcommand*{\lmat}[1]{\mathrm{#1}}
\newcommand*{\lele}[1]{\mathrm{#1}}

\theoremstyle{remark}
\newtheorem{remark}{Remark}[section]

\newcommand{\com}[1]{\iffalse{#1}\fi}

\title{A polarization tensor approximation for the Hessian in iterative solvers for non-linear inverse problems}
\author{
	\name{F.~M. Watson\textsuperscript{a,1}\footnote{\textsuperscript{1}At the time of writing F. M. Watson was at the Defence Science and Technology Laboratory. He has since moved to Thales UK Ltd.}
	, M.~G. Crabb\textsuperscript{b,2}\footnote{\textsuperscript{2}M.~G. Crabb was at the Department of Electrical and Electronic Engineering, University of Manchester. He has since moved to the School of Biomedical Engineering and Imaging Sciences, King’s College London.}, W.~R.~B. Lionheart\textsuperscript{c}, }
	\affil{\textsuperscript{a}Radar Applied Research Team, Defence Science and Technology Laboratory, Salisbury, UK;\\\textsuperscript{b}Electrical and Electronic Engineering, University of Manchester, UK; \\\textsuperscript{c}Department of Mathematics, University of Manchester, UK.}
}

\begin{document}

\maketitle
\begin{abstract}
	For many inverse parameter problems for partial differential equations in which the domain contains only well-separated objects, an asymptotic solution to the forward problem involving `polarization tensors' exists.  These are functions of the size and material contrast of inclusions, thereby describing the saturation component of the non-linearity.  As such, these asymptotic expansions can allow fast and stable reconstruction of small isolated objects.  In this paper, we show how such an asymptotic series can be applied to non-linear least-squares reconstruction problems, by deriving an approximate diagonal Hessian matrix for the data misfit term.  
	
	Often, the Hessian matrix can play a vital role in dealing with the non-linearity, generating good update directions which accelerate the solution towards a global minimum which may lie in a long curved valley, but computational cost can make direct calculation infeasible.  Since the polarization tensor approximation assumes sufficient separation between inclusions, our approximate Hessian does not account for non-linearity in the form of lack of superposition in the inverse problem.  It does however account for the non-linear saturation of the change in the data with increasing material contrast.  We therefore propose to use it as an initial Hessian for quasi-Newton schemes.

	This is demonstrated for the case of electrical impedance tomography in numerical experimentation, but could be applied to any other problem which has an equivalent asymptotic expansion.  We present numerical experimentation into the accuracy and reconstruction performance of the approximate Hessian, providing a proof of principle of the reconstruction scheme.

\end{abstract}


\section{Introduction}
\begin{figure}
	\centering
	\subfloat[Saturation of a single datum \(\lele{d}\) over material contrast \(\sigma\).\label{fig: nonlinearity material}]{\resizebox*{0.4\textwidth}{!}{
			\includegraphics[width=\textwidth]{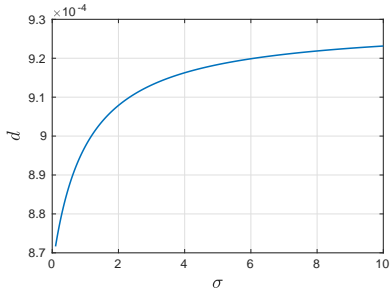}}}
	\hspace{2em}
	\subfloat[Deviation from superposition effect with separation distance between the boundaries of two inclusions.]{\resizebox*{0.4\textwidth}{!}{
			\includegraphics[width=\textwidth]{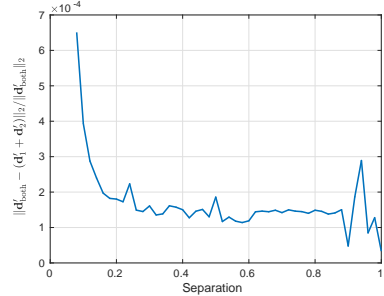}}}
	\caption{Two types of non-linearity observed in inverse boundary-value problems for PDEs, in the case of EIT. (a) Saturation of a datum \(d\) with material contrast of a single inclusion. (b) Lack of linear superposition, \(\|\lvec{d}'_{\mathrm{both}} - (\lvec{d}'_1 + \lvec{d}'_2)\|_2/\|\lvec{d}'_{\mathrm{both}}\|_2\) against the separation distance between the boundary of two inclusions. \(\lvec{d}'_{\mathrm{both}}\) is the change in data from a homogeneous domain with two inclusions, \(\lvec{d}'_1\) and \(\lvec{d}'_2\) the change in data with the first and second inclusion alone, respectively.}
	\label{fig: nonlinearity}
\end{figure}

There are two types of non-linearity in inverse parameter problems for boundary value partial differential equations (PDEs) which are observed by practitioners: lack of superposition in measurements due to multiple nearby objects; and saturation of the change in data with increasing material contrast.  These two phenomena can be observed in \figurename~\ref{fig: nonlinearity}, for the case of electrical impedance tomography (EIT).  For several PDEs in which the domain contains small well-isolated objects, a solution to the forward problem exists in the form of an asymptotic expansion involving generalised polarization tensors (GPTs), also referred to as `polarizability' tensors; see for example Ammari \textit{et al}\cite{Ammari05}.  The terms of such series are in increasing order with respect to the size of inclusion, and with GPTs that are increasingly higher rank tensors. Such asymptotic approximations exist for Maxwell's equations, metal detection, EIT, acoustics and elasticity, see for example \cite{Ledger13b, ammari2014reconstruction, LionheartKhairuddin17, Ledger19, Ammari05, Ammari07b, Kleinman86}. The polarization tensors depend on the physical regime, but in each case are a function of the shape of inclusions and material contrast. The use of such models to solve the inverse problem can then deal well with non-linearity due to material contrast.  In general they provide no information about the lack of superposition, due to assuming sufficiently separated objects.

Where the inverse problem does not only involve well isolated objects, one often attempts to solve numerically the more generic (possibly large-scale or high-dimensional) non-linear least-squares reconstruction problem
\begin{equation}
\lvec{m}_{\mathrm{im}} = \argmin_{\lvec{m}\in M} \|\mathcal{F}(\lvec{m}) - \lvec{d}\|_2^2 + \lambda R(\lvec{m}),
\label{eq: inverse prob}
\end{equation}
where \(\lvec{m}\in M\subseteq\mathds{R}^{N_m}\) is a discretisation of physical parameters to be recovered, \(\lvec{m}_{\mathrm{im}}\) the reconstructed parameters (or \textit{image}), \(\lvec{d}\in\mathds{R}^{N_d}\) the observed data, and \(\mathcal{F}:\mathds{R}^{N_m}\rightarrow\mathds{R}^{N_d}\) is the forward operator simulating data. Prior knowledge is incorporated through the regularisation term \(R:\mathds{R}^{N_m}\rightarrow\mathds{R}^+\), which stabilises ill-conditioned and over- or under-determined inverse problems, effectively preventing over-solving and fitting noise, with regularisation parameter \(\lambda>0\); see for example Tarantola\cite{Tarantola82} or Vogel\cite{Vogel02}.
Solving this problem can have a high computational cost, since iterative solution involves calculating \(\mathcal{F}(\lvec{m}^{[k]})\) for many different iterates \(\lvec{m}^{[k]}\) (the \textit{forward problem}).  Both lack of superposition and saturation must be dealt with by the reconstruction scheme.

Due to both non-linearity and ill-posedness, the topography of the cost function is often characterised by an elongated curved valley, as well as having multiple local minima\cite{Martinez12}.  The Hessian matrix of the least-squares data-misfit function can play an important role in dealing with these features of (\ref{eq: inverse prob}).  It describes parameter illumination in the data, interactions between two nearby inclusions in off-diagonal terms, and the non-linear saturation effect in the leading diagonal -- to second order.  Incorporating the Hessian (or an approximation) in a Newton-type method therefore acts to refocus the gradient direction\cite{Pratt97, Haber2000, BuiThanh12i, BuiThanh12ii, BuiThanh12iii}, often arriving more rapidly towards the minimum of the valley.  However, calculation of this large and possibly dense matrix in full is not always possible due to memory and computation time constraints.  Several alternative approaches may be taken to efficiently incorporate information from the Hessian in the update direction without calculating it directly.  These include Gauss-Newton\cite{Schweiger05, Kaltenbacher10}, quasi-Newton \cite{Haber2000, Haber05, Lavoue14}, and inexact Newton methods\cite{Santosa87, Jin12, Metivier13}, to reference a few.

In this paper, we use a polarization tensor approximation (specifically, the classical P\'{o}lya-Szeg\"{o} tensor) to derive an approximate diagonal of the Hessian matrix, which is computationally very cheap compared with calculating the true diagonal.  We propose to use this as an initial Hessian estimate in quasi-Newton schemes in a novel `mixed-model' approach.  We investigate the effectiveness of encoding information about the non-linear saturation as well as parameter illumination into the reconstruction scheme in this way.  This approach is widely applicable to any inverse parameter problem for boundary value PDEs for which such an asymptotic expansion exists.  

To demonstrate the method we apply it to the EIT reconstruction problem. This is both severely ill-posed and non-linear and so provides a sufficiently challenging test of the method.  We provide numerical evidence of the performance of the method, both in terms of the quality and efficiency of the reconstruction, as well as how accurate an approximation to the true Hessian we have.  The results presented are for 2D reconstruction problems, but it is large-scale inverse problems, where memory is a limiting factor, for which we expect most utility to be gained.  Previously, we have also used this method for 3D reconstruction of ground-penetrating radar data\cite{Watson15thesis, Watson16conf}.  Thus, we provide a \textit{proof of principle} for the method to be used for large scale non-linear inverse problems in general.

This paper is organised as follows.  In Section~\ref{sec: EIT}, we outline the theoretical background to EIT, which is the physical problem we will use to demonstrate the reconstruction method. In Section~\ref{sec: PT}, we cover the relevant literature results on asymptotic solution to the generalised Laplace equation in terms of polarization tensors, as well as properties of the classical tensor of P\'{o}lya and Szeg\"{o}. In Section~\ref{sec: recon} we discuss reconstruction schemes, and in particular Section~\ref{sec: Hess approx} uses the asymptotic approximation of Section~\ref{sec: PT} to derive the approximate diagonal Hessian matrix which we propose to use as an initial Hessian approximation in l-BFGS.  Finally, numerical experimentation and discussion of the results is presented in Section~\ref{sec: numerical}.  This includes both qualitative and quantitative comparison of the approximate Hessian to the true one, as well as reconstruction results using the approximate Hessian initialised l-BFGS scheme.\

\section{Electrical Impedance Tomography}\label{sec: EIT}
The aim of EIT is to reconstruct the conductivity of an object from low-frequency electrical measurements obtained on the boundary. In the zero frequency limit, the relationship between the electrical potential, $u$, and conductivity, $\sigma$, is governed by a second order, linear PDE
\begin{equation} 
\left\{\begin{array}{l}\nabla\cdot\left(\sigma\nabla u\right) =0\quad x\in\Omega\subset\mathds{R}^d \\
\left.\sigma \frac{\partial u}{\partial \nu}\right|_{\partial\Omega}=g  \end{array}\right.
\label{eq:EIT_common}
\end{equation}
where the domain \(\Omega\) has boundary \(\partial\Omega\) and \(d=2,3\).  In general $\sigma$ can be anisotropic but in this paper we will consider isotropic conductivity $\sigma : \Omega \to 
\mathds{R}$, bounded above and below, $0 < c < \sigma(x) < C$ almost everywhere for constants $c, C \in \mathds{R}$, denoted $\sigma \in L^{\infty}_{+}(\Omega)$. Given Neumann boundary conditions i.e. $\sigma \frac{\partial u}{\partial \nu}|_{\partial \Omega} = g$, $g \in H^{-\frac{1}{2}}(\partial \Omega)$,  \eqref{eq:EIT_common} has a unique weak solution $u \in H^{1}(\Omega) / \mathds{R}$, where the quotient space reflects that the electric potential is only defined up to a constant. We denote the Neumann-to-Dirichlet map $\Lambda_{\sigma} : H^{-\frac{1}{2}}(\partial \Omega) \to H^{\frac{1}{2}}(\partial \Omega)$ as
\begin{equation}
\Lambda_{\sigma}(g) = u|_{\partial \Omega}.
\end{equation}
The mathematical formulation of the inverse conductivity problem is then to study the determination of $\sigma$ from $\Lambda_{\sigma}$. This problem is non-linear and severely ill-posed, and in practice one only has partial knowledge of the boundary data, $\Lambda_{\sigma}$, which is subject to measurement noise (see e.g. \cite{BillChapterHolder,BorceaReview2002,UhlmannReview2009} for review articles on reconstruction algorithms and theoretical results for EIT).

\subsection{EIT forward modelling}\label{sec: EIT forward}
In practice only a finite number of currents and voltage measurements can be applied and measured respectively and on finite sized electrodes, and in many applications there is a power drop due to formation of contact impedance at electrode/domain interfaces. We consider $N_L$ electrodes and denote with  $E_l \subset \partial \Omega$ the subset of boundary in contact with the $l^{th}$ electrode and $E_{T}^{'} := \partial \Omega \setminus \bigcup_{l=1}^{N_L} E_{l}$. The complete electrode model (CEM) \cite{Somersalo:1992:EUE:132103.132112} consists of the conductivity equation (\ref{eq:EIT_common}) along with the boundary conditions
\begin{equation}
(u + \eta_{l} \sigma \frac{\partial u}{\partial \nu}) |_{E_{l}} = U_{l}, \hspace{3 mm} l=1,\ldots,L, 
\hspace{6 mm}
\sigma \frac{\partial u}{\partial \nu} = 0, \hspace{3 mm} x \in E_{T}^{'},
\hspace{6 mm}
\int_{E_{l}} \sigma \frac{\partial u}{\partial \nu} \de s = I_{l},  
\label{eq:complete_elec_bcs}
\end{equation}
where $I:=(I_1,\ldots,I_{N_L})^T \in \mathds{R}^{N_L}$, with $\sum_{l=1}^{N_L}I_{l} = 0$, are the inflow currents, $\eta:=(\eta_1,\ldots,\eta_{N_L})^T \in \mathds{R}^{N_L}$ are the contact impedances and $U:=(U_1,\ldots,U_{N_L})^T \in \mathds{R}^{N_L}$ are the potentials. The CEM forward problem is then: given $\sigma$ and $I$, determine $u$ and $U$. There is a unique weak solution $(u,U)\in (H^{1}(\Omega) \oplus \mathds{R}^{N_L}) / \mathds{R}$ to this problem \cite{Somersalo:1992:EUE:132103.132112}.

We now discuss a numerical solution to the CEM forward problem using the finite element method (FEM) (for further details on FEM see e.g. \cite{fem_silvester,fem_ciarlet,fem_zienkiewicz,VauhkonenPHD,ledger_hpfem,CrabbpFEM2017}). In the FEM the domain $\Omega$ is decomposed into $N_E$ disjoint elements $\{ \Omega_i \}_{i=1}^{N_E}$, chosen here to be triangles, with $\Omega = \cup_{i=1}^{N_E} \Omega_{i}$, joined at $N_N$ vertex nodes $\{x_i\}_{i=1}^{N_N}$. A piecewise constant conductivity discretisation, $\sigma = \sum_{i=1}^{N_E} m_i \chi_i$, is chosen here where $\chi_i$ is the characteristic function of $\Omega_i$ and ${\lvec{m}}:=(\lele{m}_1,\ldots,\lele{m}_{N_E})^T \in \mathds{R}^{N_{E}}$ is the representation of the conductivity. In the FEM, a continuous approximation to the weak solution of \eqref{eq:EIT_common}, $u_{h}$, is sought where $u \approx u_{h} = \sum_{i=1}^{N_N} \upalpha_i \psi_i$, where $\upalpha:=(\upalpha_1,\ldots,\upalpha_{N_N})^T \in \mathds{R}^{N_N}$ and $\{\psi_{i}\}_{i=1}^{N_N}$ are the shape functions, chosen here to be piecewise linear. Additionally, an approximation  $\upbeta:=(\upbeta_1,\ldots,\upbeta_{N_L})^T \in \mathds{R}^{N_L}$ to the electrode potentials is sought such that $U \approx \upbeta$. As described by Vauhkonen \cite[pg. 44]{VauhkonenPHD}, $(\upalpha,\upbeta)$ is given by the solution of the system of equations
\begin{align}
\lmat{S}(\lvec{m})
\left (
\begin {array}{c}
\lvec{\upalpha} \\
\noalign{\medskip}
\lvec{\upbeta}
\end {array}
\right)
:=
\left (
\begin {array}{cc}
\lmat{A}(\lvec{m})+\lmat{B}&\lmat{P}\\
\noalign{\medskip}
\lmat{P}^{T}&\lmat{Q}
\end {array}
\right)
\left (
\begin {array}{c}
\lvec{\upalpha} \\
\noalign{\medskip}
\lvec{\upbeta}
\end {array}
\right)
= 
\left (
\begin {array}{c}
\lmat{0}\\
\noalign{\medskip}
\lmat{I}
\end {array}
\right),
\label{eq:linear_comp1}
\end{align}
where the stiffness matrix $\lmat{A}(\lvec{m})$ is given by
\begin{equation}
\lmat{A}_{ij} := \int_{\Omega} (\sigma \nabla \psi_{i}) \cdot \nabla \psi_{j}\de\Omega, \hspace{4 mm} i,j=1,\ldots,N_N,
\label{eq:A_mat}
\end{equation}
and the matrices $\lmat{B}$, $\lmat{P}$ and $\lmat{Q}$ by
\begin{equation}
\begin{aligned}
\lmat{B}_{ij} := \sum_{l=1}^{N_L} \frac{1}{\eta_{l}}\int_{E_{l}} \psi_{i} \psi_{j} \de s, \hspace{3 mm} \lmat{P}_{jl} := -\frac{1}{\eta_{l}} \int_{E_{l}} \psi_{j} \de s, \hspace{3 mm} \lmat{Q}_{ll} := \frac{|E_{l}|}{\eta_{l}}, \\ i,j=1,\ldots,N_N, \hspace{3 mm} l=1,\ldots,N_L.
\end{aligned}
 \label{eq:CEM_BCD_mats}
\end{equation}
This can be written compactly as
\begin{equation}
\lmat{S}(\lvec{m}) \lvec{u} =  \lvec{h}
\label{eq:linear_comp_schematic}
\end{equation}
where $\lvec{u}:= (\lvec{\upalpha},\lvec{\upbeta})^{T}$, $\lvec{h} := (\lvec{0},\lvec{I})^{T}$. The potential $\lvec{u}$ is only defined up to a constant resulting in a $1$-dimensional null space of $S$, but this problem is resolved by choosing an interior node, with coordinate $x_c$, to be at zero potential $\lvec{u}(x_c)=\lvec{0}$. The $c^{th}$ row and column of $\lmat{S}$ and $c^{th}$ row of $\lvec{u}$ and $\lvec{h}$ are removed to generate $\tilde{\lmat{S}}$, $\tilde{\lvec{u}}$ and $\tilde{\lvec{h}}$, and an $N_N+N_L-1$ dimensional linear system $\tilde{\lmat{S}}(\lvec{m}) \tilde{\lvec{u}} = \tilde{\lvec{h}}$ is solved for $\tilde{\lvec{u}}$.  

The forward problem is solved with $N_L-1$ right hand side vectors $\{ \tilde{\lvec{h}}^{[n]} \}_{n=1}^{N_L-1}$, that are determined from $N_L-1$ linearly independent current vectors $\{ I^{[n]} \}_{n=1}^{N_L-1}$, yielding solutions $\tilde{\lvec{u}}^{[n]}=\tilde{\lmat{S}}(\lvec{m})^{-1}\tilde{\lvec{h}}^{[n]}$ and measured voltage patterns $U^{[n]}$, $n=1,\ldots,N_L-1$. We define the \(r\)\textsuperscript{th} measurement as the voltage difference between electrode $l$ and electrode $(l \operatorname{mod} N_L) +1$ at the application of the \(n\)\textsuperscript{th} current, $U^{[n]}_{(l \operatorname{mod} N_L)+1}-U^{[n]}_{l}$, $l=1,\ldots,N_L$, and $r=N_L(n-1)+l$. This can be written using a linear measurement operator $\lvec{t}^{[l]}$ (a column vector) that generates the \(r\)\textsuperscript{th} component of simulated data, $\lele{f}_r$, through $\lele{f}_{r} := (\lvec{t}^{[l]})^{T}\tilde{\lvec{u}}^{[n]}$, and  \(\mathcal{F}(\lvec{m}) = [\lele{f}_1,\ldots,\lele{f}_{N_d}]^T\).  This results in $N_d=\frac{1}{2}N_L(N_L-1)$ independent measurements (with a factor of $\frac{1}{2}$ accounting for redundancy in measurements due to $\Lambda_{\sigma}$ being self-adjoint, $\Lambda_{\sigma}=\Lambda_{\sigma}^{*}$).  Hereon, the tilde notation for the modified system will be dropped, and it is understood that we are using the modified system.  

\section{Generalized Polarization tensors}\label{sec: PT}
In this section, we provide some background results on generalized polarization tensors from the literature, for the expression we will make use of for our Hessian approximation. For further details, see Ammari and Kang\cite{Ammari07} and others\cite{CedioFenga98, ammari06review, LionheartKhairuddin17, ammari2014reconstruction, Ledger19, Ledger13b, Ammari07b, Ammari05}.

We wish to describe the effect on the electric potential of a single inclusion in domain \(\Omega\) via an asymptotic expansion.  Let \(B\subset\Omega\subset\mathds{R}^d\) (the inclusion) be a bounded domain containing the point \(z\).  Let the conductivity of \(B\) be \(\gamma\), where \(0<\gamma\neq1\), with the conductivity of the background equal to \(1\) so that \(\gamma\) is the ratio between conductivity of the object and conductivity of the background.  The conductivity is thus
\begin{equation}
\sigma = 1+(\gamma-1)\chi(B),
\label{eq:cond_single_inclusion}
\end{equation}
where \(\chi(B)\)  is the characteristic function of \(B\).   Denote by \(u_0\) the field in the absence of the object i.e. the solution of \eqref{eq:EIT_common} with $\sigma \equiv 1$,
and let \(u\) the perturbed field which is the solution of (\ref{eq:EIT_common}) with $\sigma$ given as in \eqref{eq:cond_single_inclusion}.  For \(B\) sufficiently far from \(\partial \Omega\), the perturbation in electric field satisfies the asymptotic formulae\cite{CedioFenga98, Ammari07}
\begin{equation}
\Lambda_D(g)(x) - \Lambda_0(g)(x) = -\epsilon^d\nabla u_0(z) \mathcal{M}\nabla_z N(x,z) + \bigo(\epsilon^{d+1})
\label{eq: bded 1st order single}
\end{equation}
as \(\epsilon\rightarrow0\), where \(\epsilon\) is the size (diameter) of the inclusion, \(\Lambda_D\) and \(\Lambda_0\) are the Neumann-to-Dirichlet maps \(\Lambda_D(g) = u|_{\partial\Omega}\) and \(\Lambda_0(g)=u_0|_{\partial\Omega}\), respectively, and 
  \(N\) is the Neumann function satisfying
\begin{equation}
\left\{\begin{array}{l} \nabla^2_x N(x,z) = -\delta(x-z)\quad\mbox{in }\Omega, \\
\left.\frac{\partial N}{\partial \nu_x}\right|_{\partial\Omega} = -\frac{1}{|\partial\Omega|},\\
\int_{\partial\Omega}N(x,z)\de s(x)=0,\end{array}\right.
\label{eq: Neumann func def}
\end{equation}
where $\delta$ is the Dirac delta distribution. The first-order polarization tensor \(\mathcal{M}\)  is the classical P\'{o}lya-Szeg\"{o} tensor associated with \(B\).  It varies with the conductivity contrast of the inclusion \(\gamma\) as well as its shape \(B\), but not with the position of the inclusion \(z\).  This tensor can be explicitly computed for disks and ellipses in the plane, as well as balls and ellipsoids in three-dimensional space.  For example, if \(B\) is an ellipse whose semi-axes of length \(a\) and \(b\) are on the \(x_1\)-- and \(x_2\)-- axis and, respectively, then its P\'{o}lya-Szeg\"{o} tensor \(\mathcal{M}\) can be written in matrix form as
\begin{equation}
\mathcal{M}=(\gamma -1)|B|\begin{bmatrix} \frac{a+b}{a+\gamma b}& 0 \\ 0 & \frac{a+b}{b+\gamma a}\end{bmatrix},
\label{eq: explicit ellipse}
\end{equation}
where \(|B|\) denotes the volume of \(B\).  Moreover, the change in tensor owing to a unitary transformation of the inclusion can also be readily computed.  If \(B'\subset\Omega\) is a rotation of the ellipse \(B\) such that it is not oriented with the coordinate axes, the first-order polarization tensors \(\mathcal{M}\) and \(\mathcal{M}'\) associated with \(B\) and \(B'\) (respectively) are related by\cite{ammari06review}
\begin{equation}
\mathcal{M} =    \lmat{R}\mathcal{M}'\lmat{R}^T,
\label{eq: tensor rotation}
\end{equation}
where \(\lmat{R}\) is the rotation matrix from the coordinate axes to the principle axes of \(B'\).

\begin{remark}\label{rm: equiv tensor}
	For any given P\'{o}lya-Szeg\"{o} tensor \(\mathcal{M}\), an elliptical (in 2D) or an ellipsoidal (in 3D) inclusion can be constructed with the same tensor.  We need only construct a tensor with the correct eigenvalues via (\ref{eq: explicit ellipse}), to which we apply a unitary transformation matrix \(R\) as in (\ref{eq: tensor rotation}) to align the eigenvectors.  Effectively, this tells us that the most this second rank tensor tells us about the effect of an inclusion on a field, is the effect of its closest fitting oriented ellipse/ellipsoid would have on that field (to first order), see e.g. \cite{Khairuddin13, ammari2014reconstruction}.
\end{remark}

Equation (\ref{eq: bded 1st order single}) can be extended to the case of a domain containing multiple inclusions \(B_1,\ldots,B_m\), sufficiently separated from both the boundary and one-another, by
\begin{equation}
\Lambda_D(g)(x) - \Lambda_0(g)(x) = -\sum_{s=1}^m\epsilon^d_s\nabla u_0(z_s) \mathcal{M}^{(s)}\nabla_z N(x,z_s) + \bigo(\epsilon^{d+1})
\label{eq: bded 1st order}
\end{equation}
as \(\epsilon\rightarrow0\), where \(\mathcal{M}^{(s)}\) and \(z_s\) are the polarization tensor and centre of the inclusion \(B_s\).  Equation (\ref{eq: bded 1st order}) is also the first term in a full asymptotic series\cite[pg. 29, thm 4.1]{ammari06review}
{\small\begin{equation}
	\Lambda_D(g)(x) - \Lambda_0(g)(x) = -\sum_{s=1}^m\sum_{|\alpha|,|\beta|=1}^d \frac{\epsilon^{|\alpha|+|\beta|+d-2}}{\alpha!\beta!}(\partial^{\alpha} u_0)(z_s)\partial_z^{\beta} N(x,z_s)\mathscr{M}_{\alpha\beta}^{(s)} + \bigo(\epsilon^{2d})
	\label{eq: bded dom exp}
	\end{equation}}
as \(\epsilon\rightarrow0\), for multi-indices \(\alpha, \beta\) and a common (largest) length scale of inclusions \(\epsilon\).  Here the P\'{o}lya-Szeg\"{o} tensors \(\mathcal{M}^{(s)}0\) are replaced by Generalised Polarization Tensors (GPTs) \(\mathscr{M}_{\alpha\beta}^{(s)}\), with the first order GPTs (for \(|\alpha|=|\beta|=1\)) being \(\mathcal{M}^{(s)}\).  The GPTs can be calculated by carrying out boundary integrals about the inclusion shapes \(B_s\) of an auxiliary field from a surrogate transmission problem, as well as other equivalent formulae\cite{Ammari05, ammari06review}. They can also be calculated for inhomogeneous inclusions, i.e. providing a single GPT which can be used to describe the field perturbed by multiple nearby (or touching) inclusions. The formulae are not needed for the purpose of this paper, but note that the components of the GPTs themselves depend only on \(\gamma\) and \(B_s\), and not on the incident field \(u_0\) or the position of the object \(z\).

The GPTs therefore provide a way to describe the change in electric field due to the shape and conductivity of a set of inclusions, separating these properties of the inclusions from the incident field.  An equivalent expression exists for the free-space problem, in which the Neumann function \(N\) is replaced by the free-space Green's function\cite{ammari06review}.  Equivalent asymptotic expansions involving GPTs for other physical modalities, with different formulae for the associated GPTs\cite{Ammari07}.

In (\ref{eq: bded dom exp}) higher order terms involve derivatives of the incident field \(u_0\) and measurement (adjoint) field \(N\).  So for incident fields which are fairly uniform these terms quickly become negligible, as well as the series converging rapidly with decreasing size of inclusions.  The same holds true for the measurement (adjoint) field, by reciprocity.  This tells us that to be sensitive to an object's shape -- beyond finding the closest fitting ellipsoid -- we must use source and adjoint fields which are non-uniform.  

The Hessian approximation which we later propose will make use of the first term in the asymptotic series, namely equation (\ref{eq: bded 1st order}), as well as explicit formulae for the P\'{o}lya-Szeg\"{o} tensor \(\mathcal{M}\).  If further explicit formulae were developed for higher order terms or differently shaped inclusions, these could readily be used also.

\subsubsection*{Neumann function computation}
For some domains an analytic solution is available for the Neumann function.  For example, for \(\Omega\) a disc of radius $\rho$ it is given by\cite[pg. 44, eq. 2.58]{Ammari07}
\begin{equation}
N(x,z) = -\frac{1}{2\pi} \Big( \ln|x-z| + \ln \Big| \frac{\rho}{|x|} x - \frac{|x|}{\rho}z \Big| + \frac{\ln \rho}{\pi} \Big).
\end{equation}
Note that \(N\) will be well-defined for the purposes of (\ref{eq: bded 1st order}), which always has \(x\in\partial\Omega\) and \(z\in\Omega\) sufficiently separated from \(\partial\Omega\) (so that \(z\neq x\)).

For more general domain shapes for which there is no analytic solution available one can instead use a numerical approximation to the Neumann function, and we now describe a finite element approximation to this.  We seek $N_N$ continuous finite element approximations to $\{ N(x,z_k) \}_{k=1}^{N_N}$, $N_{h}^{(k)} = \sum_{i=1}^{N_N} \lmat{\Upgamma}_{i}^{(k)} \psi_{i}$, where the $k\textsuperscript{th}$ approximation corresponds to a delta function source supported at node $k$, and $\lmat{\Upgamma} \in \mathds{R}^{N_N \times N_N}$ with the $k\textsuperscript{th}$ column representing the approximation to $N(x,z_k)$.  The weak formulation for \eqref{eq: Neumann func def}, in conjunction with this approximation, leads to the system of equations
\begin{equation}
\left (
\begin {array}{c}
\lmat{A}\\
\noalign{\medskip}
\lmat{F}
\end {array}
\right)
\lmat{\Upgamma}
= 
\left (
\begin {array}{c}
\lmat{G}\\
\noalign{\medskip}
\lmat{0}
\end {array}
\right),
\label{eq:Neumann_FEM}
\end{equation}
where $\lmat{A} \in \mathds{R}^{N_N \times N_N}$ is the stiffness matrix \eqref{eq:A_mat} with $\sigma=1$, and $\lmat{F} \in \mathds{R}^{1 \times N_N}$ and $\lmat{G} \in \mathds{R}^{N_N \times N_N}$ have entries
\begin{equation}
\lmat{F}_{i} = \int_{\partial \Omega} \psi_i \de s, \hspace{4 mm} \lmat{G}_{ij} =  \delta_{ij} - \int_{\partial \Omega} \frac{1}{|\partial \Omega|}\psi_i \de s \hspace{4 mm} i,j=1,\ldots,N_N.
\end{equation} 
where $\delta_{ij}$ is the Kronecker delta. The system \eqref{eq:Neumann_FEM} has full column rank, and we can compute the solution through the normal equations.  We note that in two dimensions a delta function has regularity $\delta \in H^{-1-\epsilon}(\Omega)$ for all $\epsilon >0$, and the resulting Neumann function $N(\cdot,z)$, for a source supported at point $z$, will have regularity $N(\cdot,z) \in H^{1-\epsilon}(\Omega)$ for all $\epsilon > 0$ falling just short of $H^{1}$-regularity required to guarantee  of FE approximation.  However, by elliptic regularity, $N(\cdot,z)$ will be smooth ($C^{\infty}$) in complement of any open set containing $z$, and the FE approximation will converge away from each source \cite{CrabbpFEM2017}.

For the numerical results presented in Section~\ref{sec: numerical} we use the above analytic Neumann function.  For reconstruction problems in domain shapes which do not have an analytic Neumann function, one could numerically compute and store the Neumann function once.

\section{The reconstruction scheme}\label{sec: recon}
In this section, we provide a brief overview of the numerical solution of the inverse problem with Newton-type methods, including calculation of the gradient and Hessian via the adjoint state method. We then derive our approximate diagonal Hessian matrix, making use of the asymptotic approximation of Section~\ref{sec: PT}, and describe how it can be used in a quasi-Newton method as a part of our ``mixed model'' approach.

We consider only the least-squares data misfit part of the objective function
\begin{equation}
\mathcal{J}(\lvec{m}) := \frac{1}{2}\|\mathcal{F}(\lvec{m})-\lvec{d}\|_2^2,
\label{eq: misfit}
\end{equation}
and recall from Section~\ref{sec: EIT forward} that each component \(\lele{f}_r\) of simulated data \(\mathcal{F}(\lvec{m})=[\lele{f}_1,\ldots,\lele{f}_{N_d}]^T\) involves simulating a datum by solving (\ref{eq:linear_comp_schematic}) for different boundary conditions.  Ignoring for now the omitted regularisation term, we wish to minimise (\ref{eq: misfit}) using a Newton-type method, which has an update direction given by
\begin{equation}
\lmat{B}^{[k]}\lvec{p}^{[k]} = -\nabla \mathcal{J}(\lvec{m}^{[k]}) =: -\lvec{g}^{[k]},
\end{equation}
for \(\lmat{B}^{[k]}\) some approximation to the true Hessian \(\lmat{H}^{[k]} := \nabla_{\lvec{m}}^2\mathcal{J}(\lvec{m}^{[k]})\),  and gradient \(\lvec{g}^{[k]} := \nabla_{\lvec{m}}\mathcal{J}(\lvec{m}^{[k]})\).  Dropping the superscript-\([k]\) notation, formally computing the components of the gradient and Hessian yields
\begin{equation}
\lele{g}_j= \sum_{r=1}^{N_d} \lele{J}_{rj}\left(\lele{f}_r - \lele{d}_r\right),
\end{equation}
and
\begin{equation}
\lele{H}_{ij} = \sum_{r=1}^{N_d}  \left\{ \frac{\partial \lele{f}_r}{\partial \lele{m}_i}\frac{\partial \lele{f}_r}{\partial \lele{m}_j} + \frac{\partial^2 \lele{f}_r}{\partial \lele{m}_i\partial \lele{m}_j}\left(\lele{f}_r - \lele{d}_r\right)\right\},
\label{eq: true Hessian}
\end{equation}
respectively, where \(\lele{J}_{rj} = \frac{\partial \lele{f}_{r}}{\partial \lele{m}_{j}}\) are elements of the Jacobian matrix \(\lmat{J}\), \(\lele{f}_r\) the \(r\)\textsuperscript{th} component of simulated data, \(\mathcal{F}(\lvec{m}) = [\lele{f}_1,\ldots,\lele{f}_{N_d}]^T\), and \(\lele{d}_r\) the \(r\)\textsuperscript{th} component of \(\lvec{d}\).

\subsection{Classical adjoint field formulation}\label{sec: Hess approx}
In solving the inverse problem, these derivatives are often calculated via an adjoint field formulation (see e.g. \cite{Pratt97, Plessix06}). In the continuous setting the EIT forward problem, the map $\mathrm{F}: L^{\infty}(\Omega) \to \mathcal{L}(H^{-\frac{1}{2}}(\partial \Omega),H^{\frac{1}{2}}(\partial \Omega)), \hspace{1mm} \sigma \mapsto \Lambda_{\sigma}$, is Fr\'{e}chet differentiable with respect to $L^{\infty}$ conductivity perturbations up to arbitrary order (see e.g. \cite{Dobson_Hessian_EIT,Schotland_Born_EIT}). The 1st and 2nd Fr\'{e}chet derivatives of $\mathrm{F}$ at $\sigma \in L^{\infty}(\Omega)$ in the directions $h_1, h_2 \in L^{\infty}(\Omega)$ are given by
\begin{align}
&D\mathrm{F}({\sigma}): L^{\infty}(\Omega) \to \mathcal{L}(H^{-\frac{1}{2}}(\partial \Omega),H^{\frac{1}{2}}(\partial \Omega)), \hspace{1mm} h_1 \mapsto D\Lambda_{\sigma}[h_1], \\
&D^2\mathrm{F}({\sigma}):  L^{\infty}(\Omega) \times L^{\infty}(\Omega) \to \mathcal{L}(H^{-\frac{1}{2}}(\partial \Omega),H^{\frac{1}{2}}(\partial \Omega)), \hspace{1mm} (h_1,h_2) \mapsto D^2\Lambda_{\sigma}[h_1,h_2], 
\end{align}
where $D\Lambda_{\sigma}[h_1](g) = \delta u|_{\partial \Omega}$ and $D^2\Lambda_{\sigma}[h_1,h_2](g) = \delta^2 u|_{\partial \Omega}$, with $\delta u$ and $\delta^2 u$ given by
\begin{align}
&\nabla\cdot(\sigma\nabla \delta u) = -\nabla \cdot ( h_1 \nabla u), \hspace{4 mm} &&\left.\sigma \frac{\partial \delta u}{\partial \nu}\right|_{\partial\Omega}=0, \label{eq:EIT_interior_deriv} \\
&\nabla\cdot(\sigma\nabla \delta^2 u) = -\nabla \cdot ( h_2 \nabla \delta u), \hspace{4 mm} &&\left.\sigma \frac{\partial \delta^2 u}{\partial \nu}\right|_{\partial\Omega}=0, \label{eq:EIT_interior_2ndderiv}
\end{align}
$u$ is the solution to the forward problem \eqref{eq:EIT_common}, and by definition 
{\small\begin{align}
&\int_{\partial \Omega} g^*(\Lambda_{\sigma+h_1}(g) - \Lambda_{\sigma}(g))  \de s = \int_{\partial \Omega} g^* D\Lambda_{\sigma}[h_1](g) \de s +  o(||h_1||_{L^{\infty}}), \\
&\int_{\partial \Omega} g^*(D\Lambda_{\sigma+h_2}[h_1](g)) - D\Lambda_{\sigma}[h_1](g))  \de s = \int_{\partial \Omega} g^*D^2\Lambda_{\sigma}[h_1,h_2](g) \de s +  o(||h_2||_{L^{\infty}}).
\end{align}}
Further, efficient, adjoint-field formulae for the 1st and 2nd Fr\'{e}chet derivatives at $\sigma$ are given by
\begin{align}
&\int_{\partial \Omega} g^*(D\Lambda_{\sigma}[h_1](g) ) \de s= -\int_{\Omega} h_1 \nabla u^* \cdot \nabla u \de\Omega, 
\label{eq:EIT_frechet1} \\
&\int_{\partial \Omega} g^*(D^2\Lambda_{\sigma}[h_1,h_2](g)) \de s= - \int_{\Omega} h_2 \nabla u^* \cdot \nabla \delta u \de\Omega,
\label{eq:EIT_frechet2}
\end{align}
where $u^*$ is the solution to the adjoint problem
\begin{align} 
&\nabla\cdot\left(\sigma\nabla u^*\right) =0, \hspace{4 mm} &&\left.\sigma \frac{\partial u^*}{\partial \nu}\right|_{\partial\Omega}=g^*. \label{eq:EIT_adjoint}
\end{align}
We note $u,u^{*},\delta u \in H^{1}(\Omega)$ since these are solutions of second order, linear, elliptic PDE, and thus their partial derivatives are in $L^{2}(\Omega)$. Hence, for directions $h_1,h_2 \in L^{\infty}(\Omega)$, the formulae for the Fr\'{e}chet derivatives in \eqref{eq:EIT_frechet1} and \eqref{eq:EIT_frechet2} are well-defined.

To  compute these derivatives numerically we consider the discretised forward problem from section \ref{sec: EIT} as $\lmat{S}(\lvec{m})\lvec{u}=\lvec{h}$, and measurement  $\lele{f}_r= (\lvec{t}^{[l]})^{T} \lvec{u}^{[n]}$, $n=1,\ldots,N_L-1$, $l=1,\ldots,N_L$ and $r=N_L(n-1)+l$. Computing the first partial derivative of the discretised forward problem with respect to $\lele{m}_{i}$ yields
\begin{equation}
\lmat{S}\frac{\partial \lvec{u}}{\partial \lele{m}_i} +\frac{\partial \lmat{S}}{\partial \lele{m}_i}\lvec{u}= \lvec{0}\quad\Leftrightarrow\quad \frac{\partial \lvec{u}}{\partial \lele{m}_i} = -\lmat{S}^{-1}\frac{\partial\lmat{S}}{\partial\lele{m}_i}\lvec{u},
\end{equation}
and using the linear measurement equation  $\lele{f}_r= (\lvec{t}^{[l]})^{T} \lvec{u}^{[n]}$, as well as \(\lvec{u}=\lmat{S}^{-1}\lvec{h}\), we have
\begin{equation}
\frac{\partial \lele{f}_r}{\partial \lele{m}_i}  = - (\lvec{t}^{[l]})^{T} \lmat{S}^{-1}\frac{\partial \lmat{S}}{\partial \lele{m}_i}\lmat{S}^{-1}\lvec{h}^{[n]}.
\label{eq: adjoint grad}
\end{equation}
Equation (\ref{eq: adjoint grad}) is the adjoint field formulation of the \(r\)\textsuperscript{th} measurement component of the Jacobian. Continuing by computing the second partial derivative of the discretised forward problem yields
\begin{equation}
\frac{\partial \lmat{S}}{\partial \lele{m}_j} \frac{\partial \lvec{u}}{\partial \lele{m}_i} + \lmat{S} \frac{\partial^2 \lvec{u}}{\partial \lele{m}_i\partial \lele{m}_j} + \frac{\partial \lmat{S}}{\partial \lele{m}_i} \frac{\partial \lvec{u}}{\partial \lele{m}_j}  + \frac{\partial^2 \lmat{S}}{\partial \lele{m}_i\partial \lele{m}_j}\lvec{u} = \lvec{0}.
\end{equation}
The last term is zero since the conductivity is piecewise constant on elements, and using the linear measurement equation $\lele{f}_r= (\lvec{t}^{[l]})^{T} \lvec{u}^{[n]}$ yields
\begin{equation}
\frac{\partial^2 \lele{f}_r}{\partial \lele{m}_i \partial \lele{m}_j} = - (\lvec{t}^{[l]})^{T} \lmat{S}^{-1}(\frac{\partial \lmat{S}}{\partial \lele{m}_j} \frac{\partial \lvec{u}^{[n]}}{\partial \lele{m}_i} +  \frac{\partial S}{\partial \lele{m}_i} \frac{\partial \lvec{u}^{[n]}}{\partial \lele{m}_j} ).
\label{eq: adjoint Hess}
\end{equation}
While calculating the gradient (\ref{eq: adjoint grad}) costs only the same as solving an additional forward problem, calculating each element of the Hessian in this way is still prohibitively expensive.  So, we proceed to derive a computationally cheap approximation.

\subsection{Polarization tensor Hessian approximation}\label{sec: PT Hess approximation}
We wish to derive an efficient expression to approximate these derivatives.   Consider a single mesh element \(\Omega_i\), with conductivity \(m_i\), away from the boundary \(\partial\Omega\), and assume we begin the reconstruction from a homogeneous background \(\lvec{m}_0\).  Let us first consider this element \(\Omega_i\) to be a single small inclusion in the domain \(\Omega\), and so it has an associated P\'{o}lya-Szeg\"{o} tensor \(\mathcal{M}^{(i)}\) which can be used to describe the change in field to first order via (\ref{eq: bded 1st order}).  \(\mathcal{M}^{(i)}\) is a function of \(\gamma=\lele{m}_i\), and the scale factor \(\epsilon^d\) in (\ref{eq: bded 1st order}) is (approximately) the area of \(\Omega_i\) in 2D, or the volume in 3D.

We propose to differentiate the asymptotic approximation (\ref{eq: bded 1st order}) to obtain a computationally cheap approximation to higher order derivatives. Noting that the derivatives of the fields \(u\) also satisfy a generalised Laplace problem, with a scaled version of \(u\) as the source term (i.e. the continuous analogue of (\ref{eq: adjoint grad})), there will also be an equivalent expression to (\ref{eq: bded dom exp}) describing the effect the inclusion has on the derivatives of \(u\).
Therefore, formally differentiating (\ref{eq: bded 1st order}), we have that
\begin{subequations}
\begin{align}
\frac{\partial \lele{f}_r}{\partial \lele{m}_i} \approx \lele{C}_{ri}:=& -\epsilon^d_i\nabla u_0(z_i)\frac{\partial\mathcal{M}^{(i)}}{\partial \lele{m}_i}\nabla_z N(x,z_i),\\
\frac{\partial^2 f_r}{\partial {\lele{m}_i}^2}\approx \lele{D}_{ri} :=&-\epsilon^d_i\nabla u_0(z_i)\frac{\partial^2\mathcal{M}^{(i)}}{\partial {\lele{m}_i}^2}\nabla_z N(x,z_i),\label{eq: second deriv}
\end{align}\label{eq: objective derivs}
\end{subequations}
with \(z_i\) the centre of element \(\Omega_i\).  Note that neither \(u_0\) nor \(N\) are functions of \(\lele{m}_i\), as they are defined as the solutions to homogeneous problems, so only derivatives of \(\mathcal{M}^{(i)}\) appear in (\ref{eq: objective derivs}).  This contrasts with the true Fr\'{e}chet derivative (\ref{eq:EIT_frechet1}), in which \(u\) and \(u^*\) do vary with \(\lele{m}_i\).
\begin{remark}\label{rm: deriv valid}
These expressions are valid for perturbed elements of the conductivity vector \(\lele{m}_i\) which are sufficiently isolated from other perturbations of the conductivity, as well as all elements when \(\lele{m}_i=1\) everywhere (or some other choice of constant background conductivity).
\end{remark}

By remark~\ref{rm: equiv tensor}, we know there is an equivalent ellipse/ellipsoid to \(\Omega_i\) with the same P\'{o}lya-Szeg\"{o} tensor.  We therefore propose to differentiate the analytic expression for the \(\mathcal{M}^{(i)}\) using the equivalent ellipse/ellipsoid to \(\Omega_i\).  For \(d=2\), with an ellipse whose semi-major and -minor axis align with the coordinate axis,  using the formula (\ref{eq: explicit ellipse}) this is given by
\begin{subequations}
\begin{align}
\left.\frac{\partial\mathcal{M}^{(i)}}{\partial \lele{m}_i}\right|_{\gamma}=& |B|\begin{bmatrix}\frac{a+b}{a+\gamma b} & 0 \\ 0 & \frac{a+b}{b + \gamma a} \end{bmatrix} + (\gamma-1)|B|\begin{bmatrix}-\frac{b(a+b)}{(a+\gamma b)^2} & 0 \\ 0 & -\frac{a(a+b)}{(b+\gamma a)^2} \end{bmatrix}, \label{eq: first deriv polya}\\
\left.\frac{\partial^2\mathcal{M}^{(i)}}{\partial {\lele{m}_i}^2}\right|_\gamma=&2|B|\begin{bmatrix}-\frac{b(a+b)}{(a+\gamma b)^2} & 0 \\ 0 & -\frac{a(a+b)}{(b+\gamma a)^2}\end{bmatrix}  + (\gamma-1)|B|\begin{bmatrix}\frac{2b^2(a+b)}{(a+\gamma b)^3} & 0 \\ 0 & \frac{2a^2(a+b)}{(b+\gamma a)^3}\end{bmatrix}. \label{eq: second deriv polya}
\end{align}
\label{eq: derivs polya}
\end{subequations}
For oriented ellipses rotated by \(\theta\) from the coordinate axis, the differentiated tensors in (\ref{eq: derivs polya})  rotate in the same way as in (\ref{eq: tensor rotation}) using the standard linear algebra rotation matrices.

\begin{figure}
	\centering
	\subfloat[First derivative of tensor components.]{\resizebox*{0.4\textwidth}{!}{
			\includegraphics[width=\textwidth]{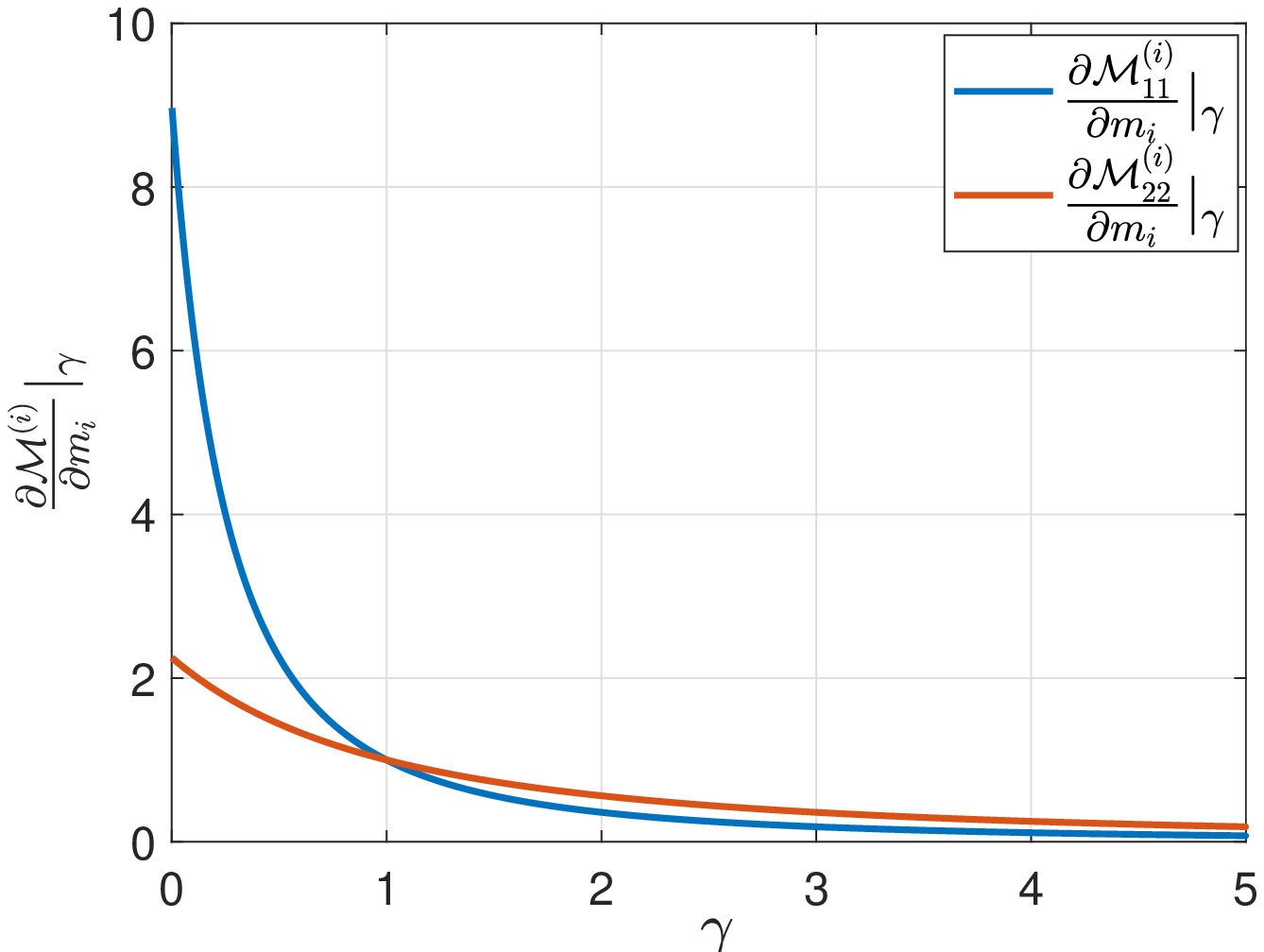}}}
		\hspace{2em}
	\subfloat[Second derivative of tensor components.]{\resizebox*{0.4\textwidth}{!}{
			\includegraphics[width=\textwidth]{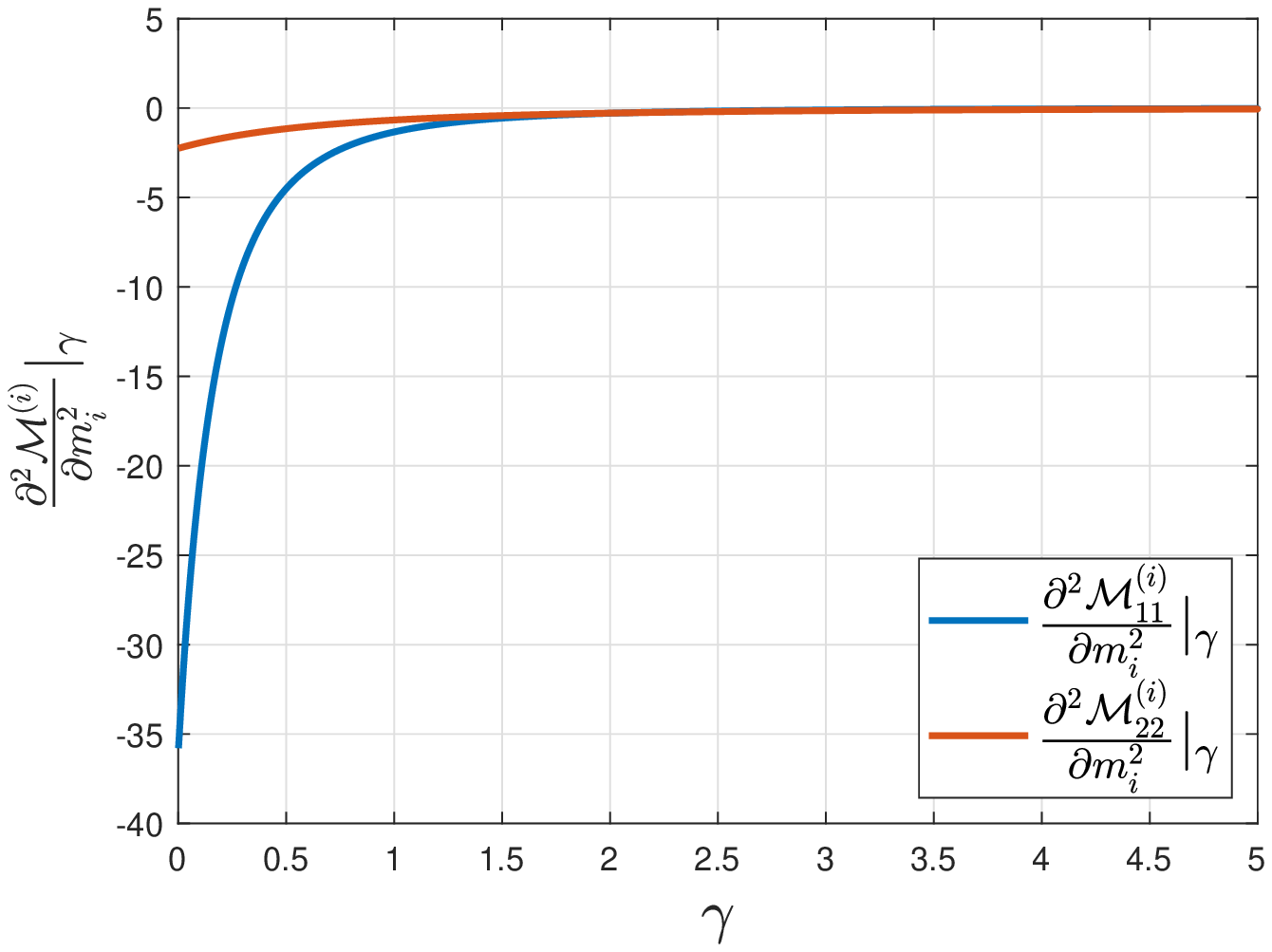}}}
		\caption{Derivative of components of the polarization tensor for an ellipsoid \(B\) with \(a=1\), \(b=2\), \(|B|=1\).}
		\label{fig: dM curves}
\end{figure}

Through the expressions in (\ref{eq: derivs polya}), we see that (\ref{eq: objective derivs}) includes the effect of saturation with material contrast depicted in \figurename~\ref{fig: nonlinearity material}: specifically, that the derivatives of tensor components tend to zero as \(\gamma\rightarrow+\infty\), and to some finite constant for \(\gamma\rightarrow0\) (which is positive for the first derivative, and negative for the second).  These components are shown in \figurename~\ref{fig: dM curves} for \(a=1, b=2\) and \(|B|=1\), and the non-linear saturation effect is clearly observed.

Recalling that ``closest fitting ellipsoid'' is electrically in the sense of (\ref{eq: explicit ellipse}), not geometrically, we note that we do not have a simple rule defining what this is for a given triangular or tetrahedral element; this is a subject of current research.  For numerical implementation we therefore must make some (possibly heuristic) choice, and for simplicity in our numerical experimentation we have chosen the Steiner inellipse\cite{SteinerWolfram},\cite[pg. 11, thm. 4.2]{Marden66}.  It is possible to calculate the exact Fr\'{e}chet derivative of GPTs for an arbitrary shaped element, see for example\cite{ammari2014reconstruction}, although this does not provide a computationally cheap tool.

The asymptotic approximations to derivatives (\ref{eq: objective derivs}) can be used to form a diagonal approximation \(\tilde{\lmat{H}}\) to the Hessian matrix,
\begin{equation}
\tilde{\lele{H}}_{ij}(\lvec{m}) = \left\{\begin{array}{ll} \sum_{r=1}^{N_d}\left\{\lele{C}_{ri}(\lele{m}_i)^2 + \lele{D}_{ri}(\lele{m}_i)(\lele{f}_r - \lele{d}_r)\right\}  \qquad & \mbox{if }i=j \\ 0 \qquad & \mbox{otherwise}.
\end{array}\right.
\label{eq: approx Hess}
\end{equation}
From remark~\ref{rm: deriv valid} above, this diagonal matrix will be a valid approximation for \({\lvec{m}-\lvec{m}_0}\) sufficiently sparse with sufficiently separated elements, where \(\lvec{m}_0\) is an initial homogeneous conductivity.  Where this is not the case, those elements corresponding to nearby perturbations will neglect non-linear interactions between the elements, in which case the approximation (\ref{eq: bded 1st order}) which this is based upon is no longer \(\bigo(\epsilon^{d+1})\) accurate.  (\ref{eq: approx Hess}) may still be a useful approximation, depending on how it is used, where these non-linear corrections are small -- which is to be determined through numerical experimentation.

The components \(\lele{C}_{ri}\) and \(\lele{D}_{ri}\) are straightforward to calculate via (\ref{eq: objective derivs}) and (\ref{eq: derivs polya}), requiring only a small number of 2- or 3-d matrix-vector multiplications each.  Moreover, the terms \(\nabla u_0\) will already have been calculated during the calculation of the \(u_r\), and, depending on the domain, an analytic expression for \(\nabla_z N(x,z)\) may be available.  This compares favourably to the full calculation (\ref{eq: adjoint Hess}), in which each element of the matrix requires another solution to the forward problem.

Equations (\ref{eq: objective derivs}) expresses the Fr\'{e}chet derivative of the EIT forward problem as the derivatives of the asymptotic series. Similarly, Ammari \textit{et al}\cite[thm. 5.1]{ammari2014reconstruction} provide an expression for the first Fr\'{e}chet derivative (\ref{eq:EIT_frechet1}) of inhomogeneous GPTs with respect to an internal conductivity of the inclusion, which is given as the difference between two asymptotic series.  This is used for a sensitivity analysis of GPTs to conductivity, to reconstruct inhomogeneous GPTs, and then subsequently reconstruct the inhomogeneous conductivity of an inclusion.  The approach we present (and its application) is subtly different -- for calculating derivatives, one could view our method as truncating the difference of these two series in \cite[thm. 5.1]{ammari2014reconstruction} to first order, and identifying the difference between these remaining series in terms of Fr\'{e}chet derivatives of polarization tensors themselves when taking the same limit.  We then propose to reconstruct an inhomogeneous domain, rather than inhomogeneous inclusions.

\begin{remark}\label{rem: mixed terms}Note that in (\ref{eq: second deriv}) it is not possible to calculate mixed derivatives, due to the summation term in (\ref{eq: bded 1st order}).  We could consider mixed derivatives for closely-spaced inclusions, using the polarization tensor for two closely spaced inclusions.  In other words, we might approximate the mixed derivative by
\begin{equation}
\frac{\partial^2 \lele{f}_r}{\partial \lele{m}_i \partial \lele{m}_j} = -\epsilon_{(ij)}\nabla U(z_{(ij)})\frac{\partial^2\mathcal{M}^{(ij)}}{\partial {\lele{m}_{(ij)}}^2} \nabla_z N(x, z_{(ij)}) + \bigo(\epsilon^{d+1}),
\end{equation}
where \(\mathcal{M}^{(ij)}\) is the polarization tensor of inclusion \(B_{(ij)} = B_i\cup B_j\), for closely spaced (joining) elements \(B_i\) and \(B_j\), with weighted centre point \(z_{(ij)}\), weighted average conductivity \(\lele{m}_{(ij)}\), and size \(\epsilon_{(ij)}\).  In order for this to provide any computational benefit, one would need an easy way to determine (approximately) the equivalent ellipse/ellipsoid of two adjoining elements to write \(\mathcal{M}^{(ij)}\) in the analytic form of the P\'{o}lya-Szeg\"{o} tensor.  Such composite object tensors have been calculated for example for the Eddy-current problem\cite{Ledger19}, but analytic expressions are not currently available.
\end{remark}

\subsection{Quasi-Newton inversion schemes}\label{sec: quasi-newton}
We consider solving the inverse problem (\ref{eq: inverse prob}) using a Newton-type method, which uses the update direction
\begin{equation}
\lmat{B}^{[k]} \lvec{p}^{[k]} = -\lvec{g}^{[k]},
\label{eq: quasi Newton}
\end{equation}
where \(\lvec{g}^{[k]}\) is the gradient for iterate \(\lvec{m}^{[k]}\) at iteration \(k\), and \(\lmat{B}^{[k]}\) some approximation to the Hessian matrix.  For a robust and efficient solution, we would expect an appropriate choice of \(\lmat{B}^{[k]}\) to have a similar structure to the true Hessian (i.e. having similar distribution of eigenvalues and eigenvectors), and require only a small number of solutions to the forward problem to calculate or update respectively. Additionally, it preferably has limited storage requirements, and the cost of solving (\ref{eq: quasi Newton}) is not so large that it outweighs the reduction in number of forward problems solved.  

We could use the diagonal approximate Hessian \(\tilde{\lmat{H}}^{[k]}\) from (\ref{eq: approx Hess}) directly in (\ref{eq: quasi Newton}), denoting \(\tilde{\lmat{H}}^{[k]}:=\tilde{\lmat{H}}(\lvec{m}^{[k]})\).  This would be computationally very cheap and ought to provide a good approximation to the spacing between contours of the objective function about the iterate \(\lvec{m}^{[k]}\), but it will not incorporate any information about the curvature of these contours.  This could result in an update direction that deals reasonably well with different parameter illumination and the ill-posedness in certain directions, but less well with the non-linearity of the inverse problem.   

We propose instead to use the approximate Hessian (\ref{eq: approx Hess}) as an initial Hessian approximation for a quasi-Newton method, updated each iteration, which is the main contribution of this paper.  Some quasi-Newton methods such as l-BFGS allow a different initial Hessian approximation to be used at each iteration\cite[pp. 177]{Nocedal00}, which also allows us to update \(\tilde{\lmat{H}}^{[k]}\) for each iterate \(\lvec{m}^{[k]}\).  We expect that such an approach should initialise the quasi-Newton method with good information about the ratio of eigenvalues at the current iterate (i.e. the spacing between contours of the objective function), allowing it to more effectively approximate the dominant eigenvectors  within fewer iterations.  This should far more rapidly build a good approximation to the true Hessian than, say, initialising with a multiple of the identity, providing information about the curvature of contours (non-linearity) which is not present in \(\tilde{\lmat{H}}^{[k]}\).  

It is worth noting that for \(\lvec{m}^{[k]}\) containing perturbations not sufficiently separated, \(\tilde{\lmat{H}}^{[k]}\) may be a poor approximation to the diagonal of the Hessian.  This may limit the effectiveness of our approach, but we could choose instead to initialise with \(\tilde{\lmat{H}}^{[0]}\) or otherwise some thresholded version.

\subsubsection{Computational cost}\label{sec: comp cost}
Recall from Section~\ref{sec: EIT forward} that \(N_E\) denotes the number of elements in the reconstruction domain, \(N_N\) the number of finite element nodes in the simulation domain, and \(N_L-1\) the number of current vectors to generate the dataset of \(N_d = \frac{1}{2}N_L(N_L-1)\) datapoints. Calculating \(\tilde{\lmat{H}}^{[k]}\) each iteration via (\ref{eq: derivs polya}), (\ref{eq: objective derivs}) and (\ref{eq: approx Hess}) involves \(\bigo(N_dN_E)\) operations as \(N_E, N_L\rightarrow\infty\).  This ignores the one-time cost of calculating the Neumann functions in the first iteration which can then be stored in memory (or permanently for repeated experiments in the same domain).  This is the same as the \(\bigo(N_d N_E)\) complexity to form the diagonal of the Gauss Newton approximate Hessian \(\diag(\lmat{J}^T\lmat{J})\), assuming \(\lmat{J}\) has already been formed and stored in calculation of the gradient of the objective function (otherwise this would be a \(\bigo(N_dN_N^3)\) complexity using an adjoint method).  

Calculating the diagonal of the true Hessian via (\ref{eq: adjoint Hess}) has computational complexity \(\bigo(N_d N_E N_N^3)\), or \(\bigo(N_d N_E^2 N_N^3)\) for the full (dense) matrix, with this high cost being driven by the need to re-simulate fields.  These complexities provide an upper bound, reducing for a sparse Hessian matrix or if an LU decomposition of the system matrix \(S\) has been calculated and stored.  The computational complexity of the true Hessian is also reduced to \(\bigo(N_L N_E N_N^3)\) when calculated via an adjoint field approach, which also has a higher memory requirement of \(\bigo(N_L N_N N_E)\).  These differing complexities and memory requirements for calculating the Hessian highlight the common trade-off between speed and memory use for iterative PDE-based reconstruction methods.  Nonetheless, it will be significantly more expensive to calculate than the proposed approximate Hessian using either method, which has no requirement to re-calculate fields.

\section{Numerical experimentation}\label{sec: numerical}
In this section, we undertake some numerical experiments into the effectiveness of the approximate Hessian to initialise l-BFGS.  We consider both the accuracy of this approximate Hessian, as well as how well the reconstruction scheme performs.  We aim to show that this provides a computationally cheap way of incorporating second derivative information into reconstruction schemes, which in some cases helps to improves contrast of images, as well as improving the rate of convergence compared to methods which use only first-derivative information.  Importantly, this method remains a feasible choice for large-scale problems in which either storing dense Hessian approximations or calculating second derivatives directly is prohibitively expensive, noting the discussion on computational cost in Section~\ref{sec: comp cost}.

Our numerical experiments were implemented in \textsc{Matlab} using the developers version of the open-source EIT package EIDORS\cite{Eidors3d,Eidors3_9}, which uses the complete electrode model outlined in Section~\ref{sec: EIT forward}.   Polarization tensor Hessian reconstructions were benchmarked against a fairly standard non-linear Gauss-Newton method to solve the optimisation problem in EIT as implemented by the EIDORS library function \texttt{inv\_solve\_core}.  The functions and scripts used in this paper can also be obtained via the EIDORS developers SVN repository. This includes the polarization tensor Hessian approximation, the Neumann functions for a disc and free-space in 2D, as well as the inversion scheme.  It also includes the calculation of the true Hessian via an adjoint method.  Experiments were carried out on a 2D disc with radius \(r_0=1\) with $N_L=16$ electrodes, of background conductivity \(\sigma_0=1\).  The meshes were created using the EIDORS function \texttt{ng\_mk\_cyl\_models}, which calls routines from Netgen Mesher\cite{netgen97}.

\subsection{Validating the Hessian approximation}

\subsubsection*{Qualitative comparison}
We begin by comparing components of the approximate Hessian terms given by equations (\ref{eq: objective derivs}), (\ref{eq: derivs polya}) and (\ref{eq: approx Hess}) to the true Hessian components (\ref{eq: true Hessian}).  For this, we take \(\Omega\) to be the 2D unit disc.  As discussed in Section~\ref{sec: PT Hess approximation}, for simplicity we use the Steiner inellipse as the closest fitting ellipse to the triangular finite elements. For comparison, we also calculate approximate Hessian matrices in which we use the free-space Green's function in place of the Neumann function -- i.e. using the equivalent expression to (\ref{eq: bded 1st order single}) for free-space.  This helps to illustrate the extent to which boundary interactions are included in the model.

\begin{figure}
	\centering
	\includegraphics[width=\textwidth]{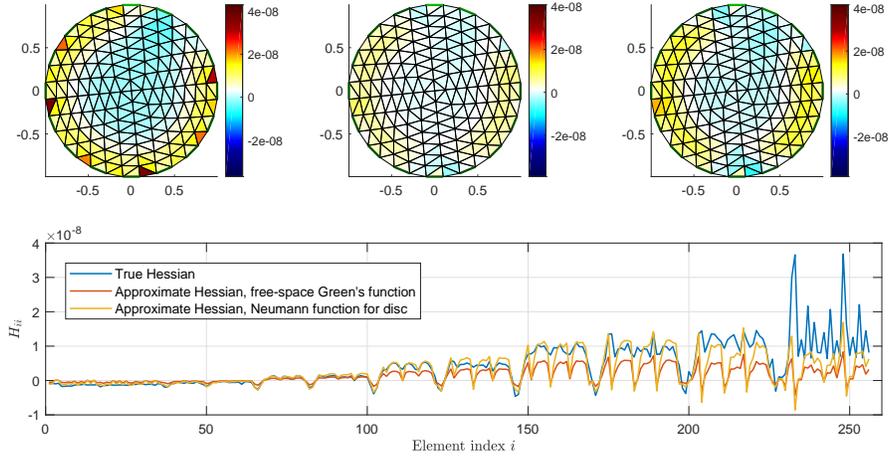}\\
	\caption{Comparison of true and approximate diagonals of the Hessian matrix calculated on a homogeneous disc, for simulated data with a single inclusion with conductivity \(\sigma=2.3\). The true Hessian is shown top-left (with elements of the matrix mapped to their corresponding element in the domain)  and in blue below.
	The approximate Hessian using the freespace Green's function is shown top-centre and in red below. The approximate Hessian using the Neumann function for the disc is shown top-right, and in yellow below.  \textit{}Higher element indices correspond to elements closer to the boundary.
	}
	\label{fig: deriv vary data}
\end{figure}

\begin{figure}
	\centering
	\subfloat[Element 26, between the inclusion and domain centre.]{\resizebox*{0.32\textwidth}{!}{\includegraphics[width=\textwidth]{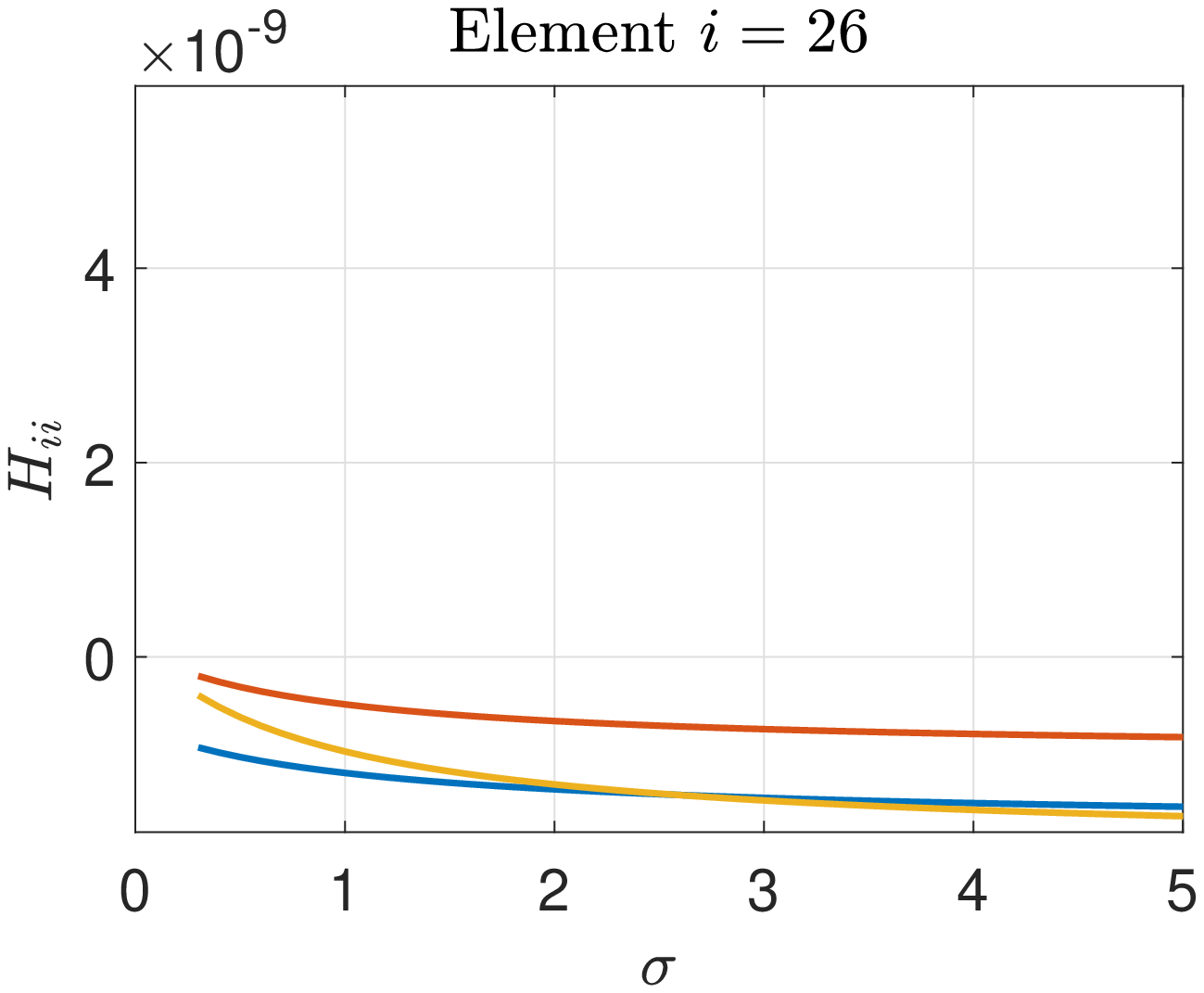}}}
	\subfloat[Element 44, within the inclusion.]{\resizebox*{0.32\textwidth}{!}{\includegraphics[width=\textwidth]{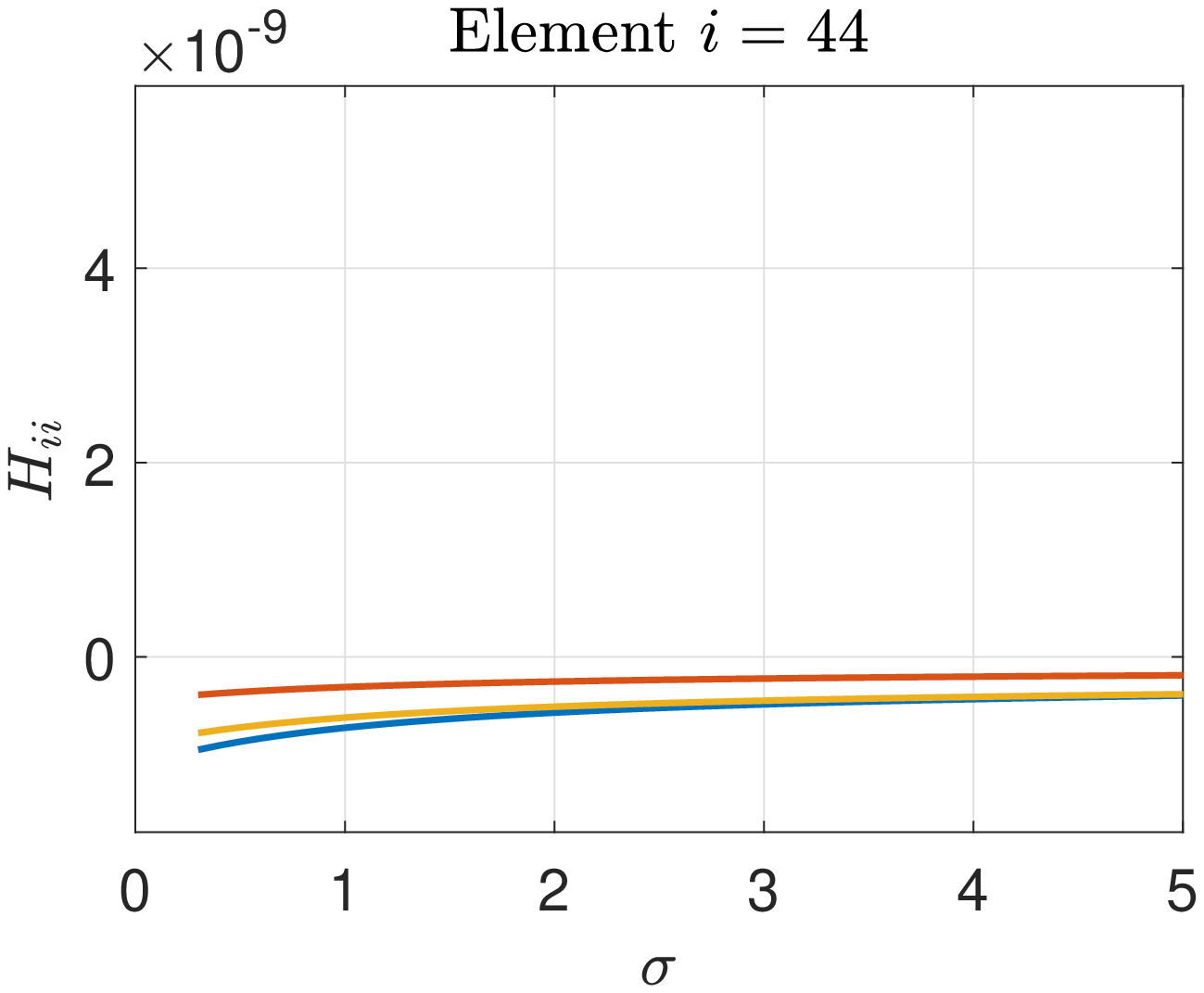}}}
	\subfloat[Element 113, between the inclusion and the boundary.]{\resizebox*{0.32\textwidth}{!}{\includegraphics[width=\textwidth]{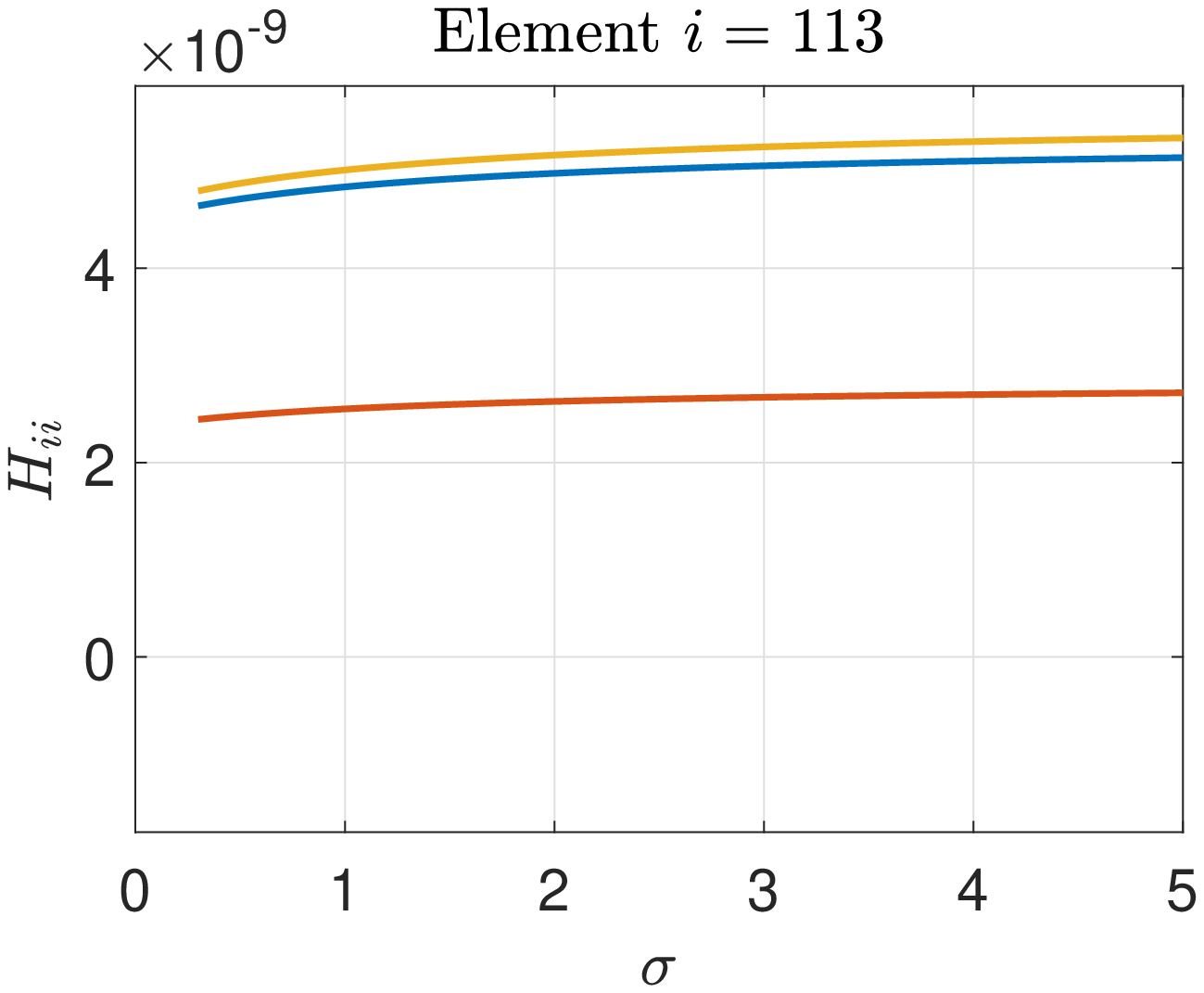}}}
	\caption{\(H_{ii}\) for selected elements \(i\) away from the boundary, as a function of inclusion conductivity \(\sigma\).  Blue shows the true Hessian, and red and yellow the approximations using the free space Green's function and Neumann function on a unit disc, respectively.}
	\label{fig: deriv saturation vary data}
\end{figure}

\figurename~\ref{fig: deriv vary data} shows the diagonal of the true and approximate Hessian matrices, calculated in a homogeneous reconstruction domain, in which data was simulated for a domain with a single circular inclusion with conductivity \(\sigma=2.3\sigma_0\), radius \(0.25\), centred at \((0.3, 0.3)^T\).  These simulations were repeated for different inclusion conductivities \(0.5\sigma_0<\sigma<5\sigma_0\), and selected elements \(H_{ii}\) (away from the boundary) are shown as a function of inclusion conductivity \(\sigma\) in Figure~\ref{fig: deriv saturation vary data}.
We see that away from the domain boundary (and in particular away from the electrodes)  the Polarization Tensor Hessian is qualitatively similar, and when using the Neumann function for a disc we have a reasonable approximation.  As the conductivity of the inclusion  \(\sigma\) in the simulation domain is varied, the Hessian approximation calculated for a homogeneous domain simply vary linearly with the data residual \(\delta\lvec{d}\) (which itself varies non-linearly with \(\sigma\)); this appears sufficient to capture the main non-linear features.

Closer to the boundary the approximation appears poor.  The indexing of elements is such that all of those from 197 to 256 share at least a node with an electrode in the FEM mesh (and so also touch the boundary), and every element from 225 has an edge on the boundary (with every other of these being an electrode).  Given this numbering, from \figurename~\ref{fig: deriv vary data} it would appear that in this case you do not need to move far from the boundary (or electrodes) before a reasonable approximation is gained.  In all cases, using the free-space Green's function (red lines in \figurename s~\ref{fig: deriv vary data}~and~\ref{fig: deriv saturation vary data}) in place of the Neumann function, and thus neglecting any effects of boundary interaction, results in a much poorer approximation apart from in the centre of the domain.

\subsubsection*{Principal angles and quantitative comparison}
To help understand the utility as well as accuracy of the approximate Hessian \(\tilde{\lmat{H}}\) initialised BFGS matrix, we look at both the relative residual between the true and BFGS Hessian in the Frobenius norm, as well as how their (dominant) singular vectors are aligned with the true matrices. The latter will suggest how closely the l-BFGS approximate Hessians will act on a gradient vector to change the update direction, closer to that of the true Hessian if the singular vectors are more closely aligned, but ignores possible scaling differences (which are easily resolved through a linesearch or within the quasi-Newton update itself).  To compare singular vectors, we use the principal angles \(\Theta\) between the subspaces they span via\cite{golub13}
\begin{equation}
\cos\Theta = \Sigma(\lmat{V}^H \lmat{W}),
\end{equation}
where \(\Sigma(\lmat{A})\) denotes the singular values of \(\lmat{A}\), and \(\lmat{V}\), \(\lmat{W}\) are the matrices whose columns are the right singular vectors.

We compare the right singular vectors associated with the 20 largest singular values.  Comparing a large number of singular vectors results in subspaces which are almost the same (since each complete set form an orthonormal basis, taking all of them would give the same subspace), whereas a small number provides little to compare and so subspaces  will generally be orthogonal.  Other choices of a similar magnitude would be equally valid.

\begin{figure}
	\centering
	\includegraphics[width=0.5\textwidth]{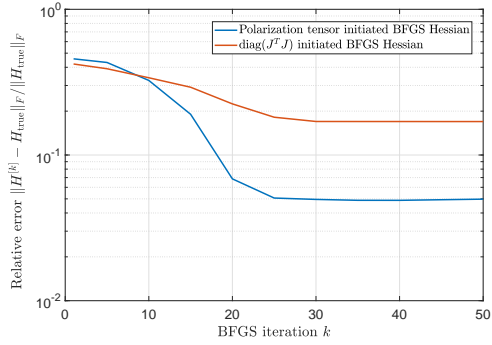}
	\caption{Relative error of the BFGS approximate Hessian to the true Hessian in Frobenius norm, initialised by the polarization tensor approximation (blue) and \(\diag(\lmat{J}^T\lmat{J})\) (red).}
	\label{fig: H convergence}
\end{figure}

\figurename~\ref{fig: H convergence} shows the relative error in Frobenius norm of the BFGS to true Hessian matrix against BFGS iteration.  The BFGS matrix initialised by \(\tilde{\lmat{H}}\) is shown in blue, and  by \(\diag(\lmat{J}^T\lmat{J})\) in red, for up to 50 BFGS iterations.  The \(\tilde{\lmat{H}}\) initialisation clearly enables a much better approximation, reaching approximately 4\% relative error after 25 iterations versus greater than 10\% for the \(\diag(\lmat{J}^T\lmat{J})\) initialisation.  In this scenario, both BFGS initialisations make very little further improvement after 25 iterations.  This is likely due to update directions \(\lvec{p}^{[k]}\) being similar across later iterations, as well as progressively smaller in norm, as the \(\lvec{x}^{[k]}\) converge to the solution, so little additional information about the Hessian (or its action on differing vectors) is gained.

\begin{figure}
	\centering
	\subfloat[\(\tilde{\lmat{H}}\) initiated BFGS Hessian.]{\resizebox*{0.48\textwidth}{!}{
		\includegraphics[width=0.98\textwidth]{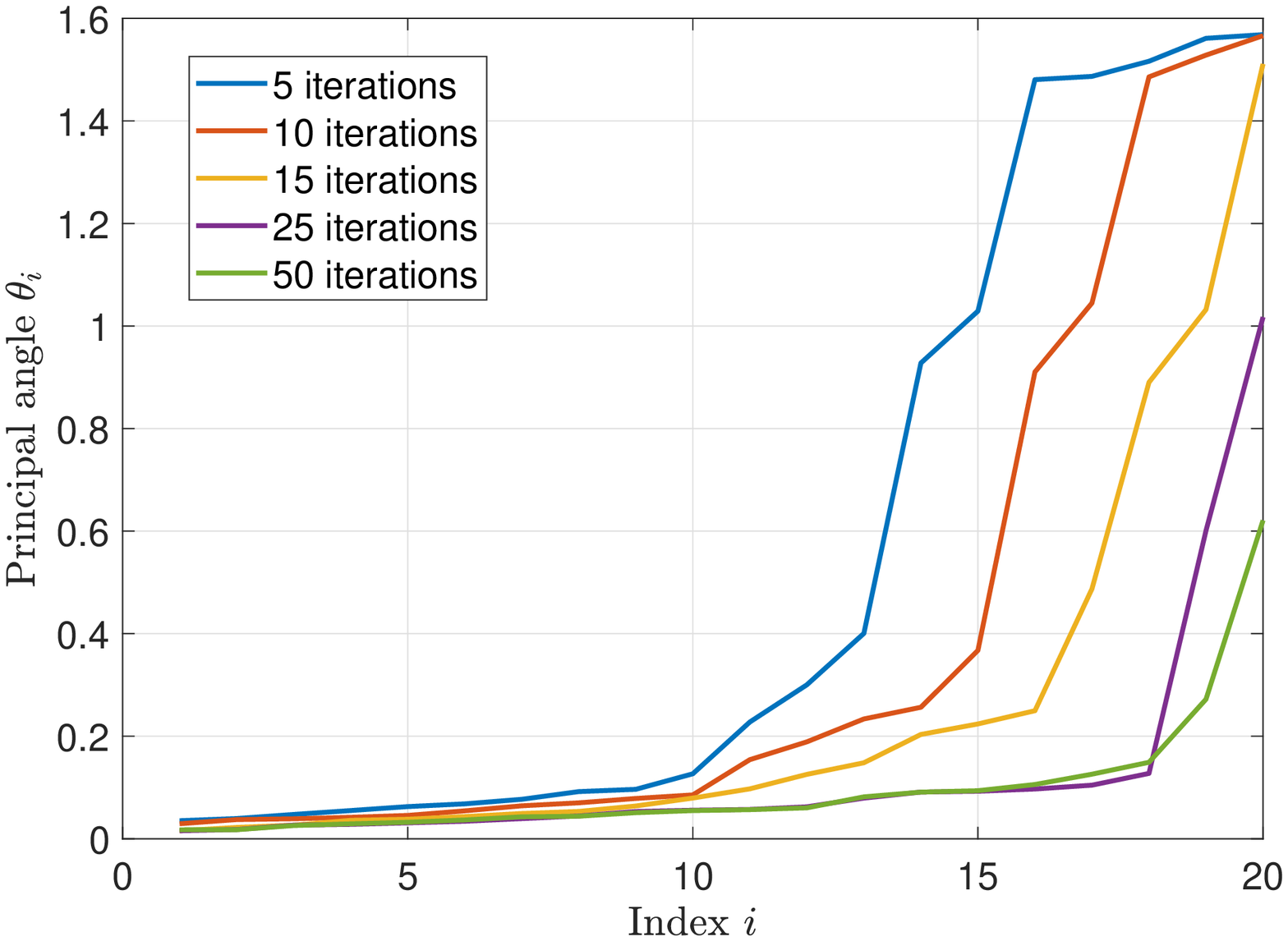}
		\label{fig: principal angles PT}}}\hspace{5pt}
	\subfloat[\(\diag(\lmat{J}^T\lmat{J})\) initiated BFGS Hessian.]{\resizebox*{0.48\textwidth}{!}{
		\includegraphics[width=0.98\textwidth]{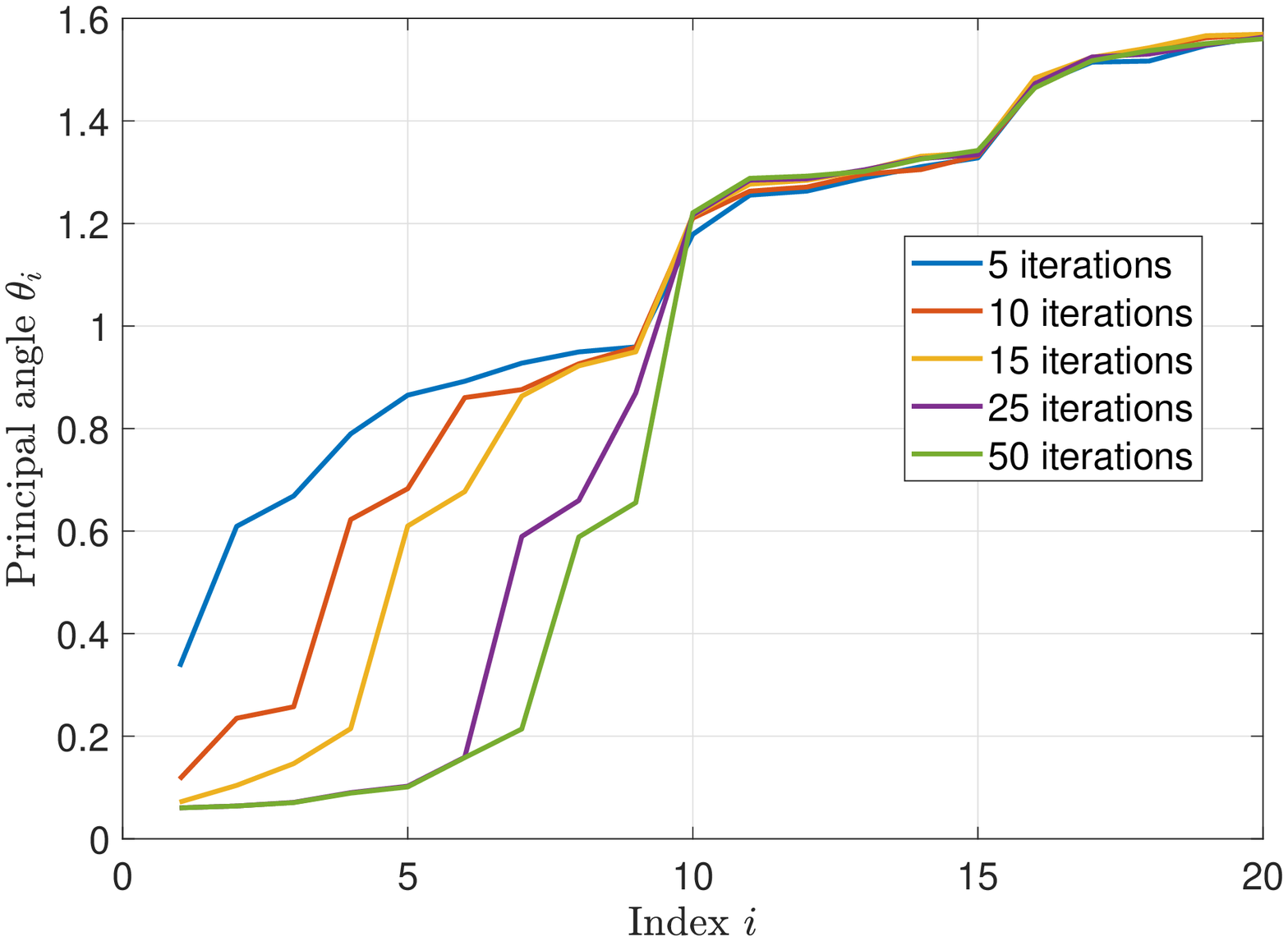}
		\label{fig: principal angles GN}}}
	\caption{Principal angles \(\theta_i\) between subspaces \(\mbox{span}(\tilde{\lmat{V}})\) and \(\mbox{span}(\tilde{\lmat{W}})\), where \(\tilde{\lmat{V}}\) and \(\tilde{\lmat{W}}\) are the matrices formed by the first 20 right singular vectors of the true Hessian \(\lmat{H}\) and BFGS approximate Hessian \(\lmat{B}^{[k]}\), respectively, after \(k\) BFGS iterations.  In \figurename~\ref{fig: principal angles PT} the BFGS approximate Hessian is initialised with the diagonal polarization tensor Hessian approximation, and in \ref{fig: principal angles GN} it is initiated with \(\diag(\lmat{J}^T\lmat{J})\).}
	\label{fig: principal angles}
\end{figure}

\figurename~\ref{fig: principal angles} shows the principal angles between the subspace spanned by the first 20 singular vectors of the true and BFGS approximate Hessians after a number of outer BFGS iterations.  In \figurename~\ref{fig: principal angles PT}, the Hessian was initialised with \(\tilde{\lmat{H}}\), and in \figurename~\ref{fig: principal angles GN} it was initialised with \(\diag(\lmat{J}^T\lmat{J})\).  It is clear from these that the \(\tilde{\lmat{H}}\) initialisation is also much better able to resolve the dominant eigendirections of the Hessian.  For example, after 5 iterations the \(\tilde{\lmat{H}}\) Hessian has 10 principal angles less than \(\pi/15\approx0.21\), but the \(\diag(\lmat{J}^T\lmat{J})\) initialisation still has only 7 similarly small angles after a full 50 iterations.

\FloatBarrier
\subsection{Reconstruction results}
\begin{figure}
	\centering
	\subfloat[True domain (left), and slice through \(x=y\) of reconstructions (right).]{\resizebox*{0.95\textwidth}{!}{
			\includegraphics[width= 0.75\textwidth]{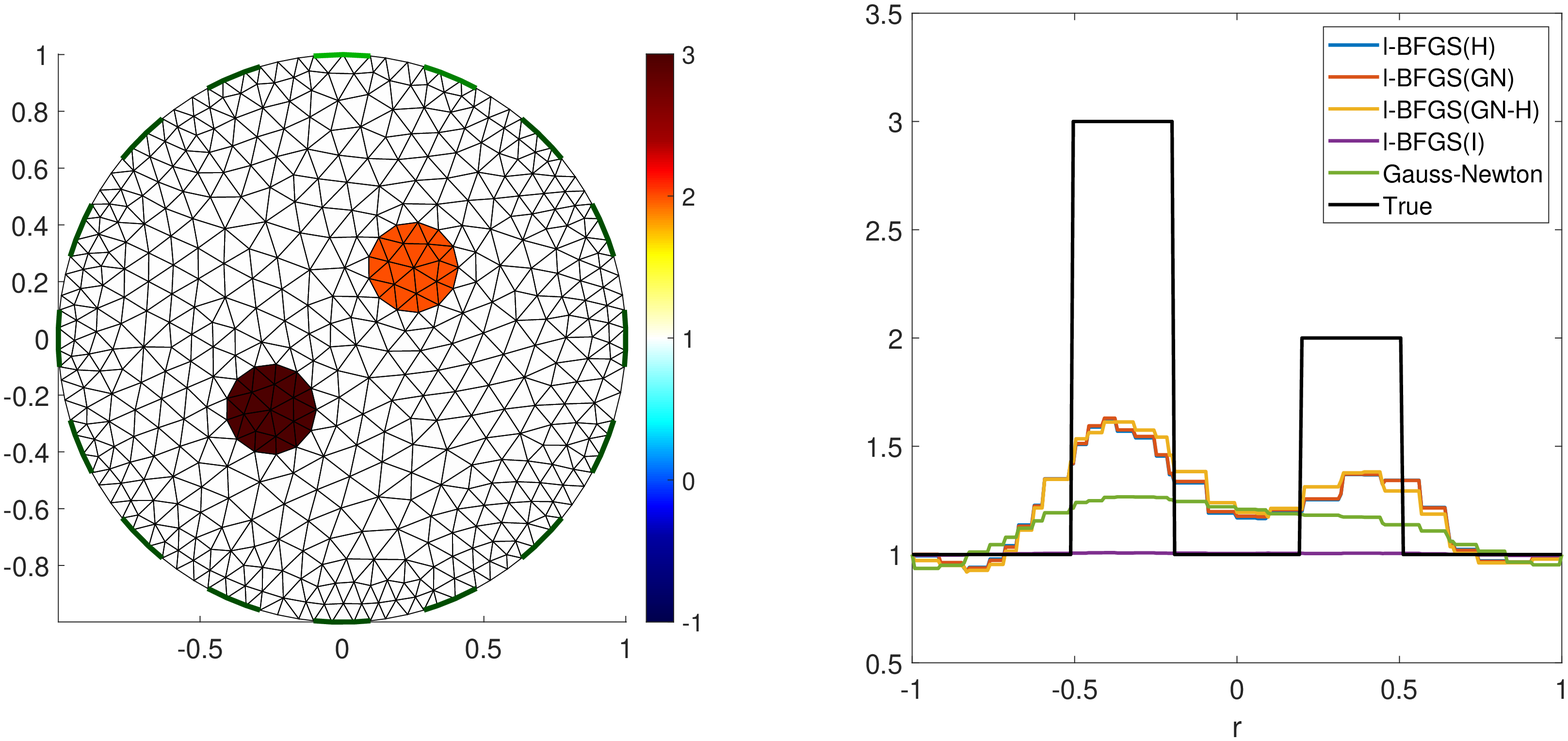}}}\\
	\subfloat[Reconstruction results with Gauss-Newton and l-BFGS]{\resizebox*{0.95\textwidth}{!}{	
			\includegraphics[width= 0.95\textwidth]{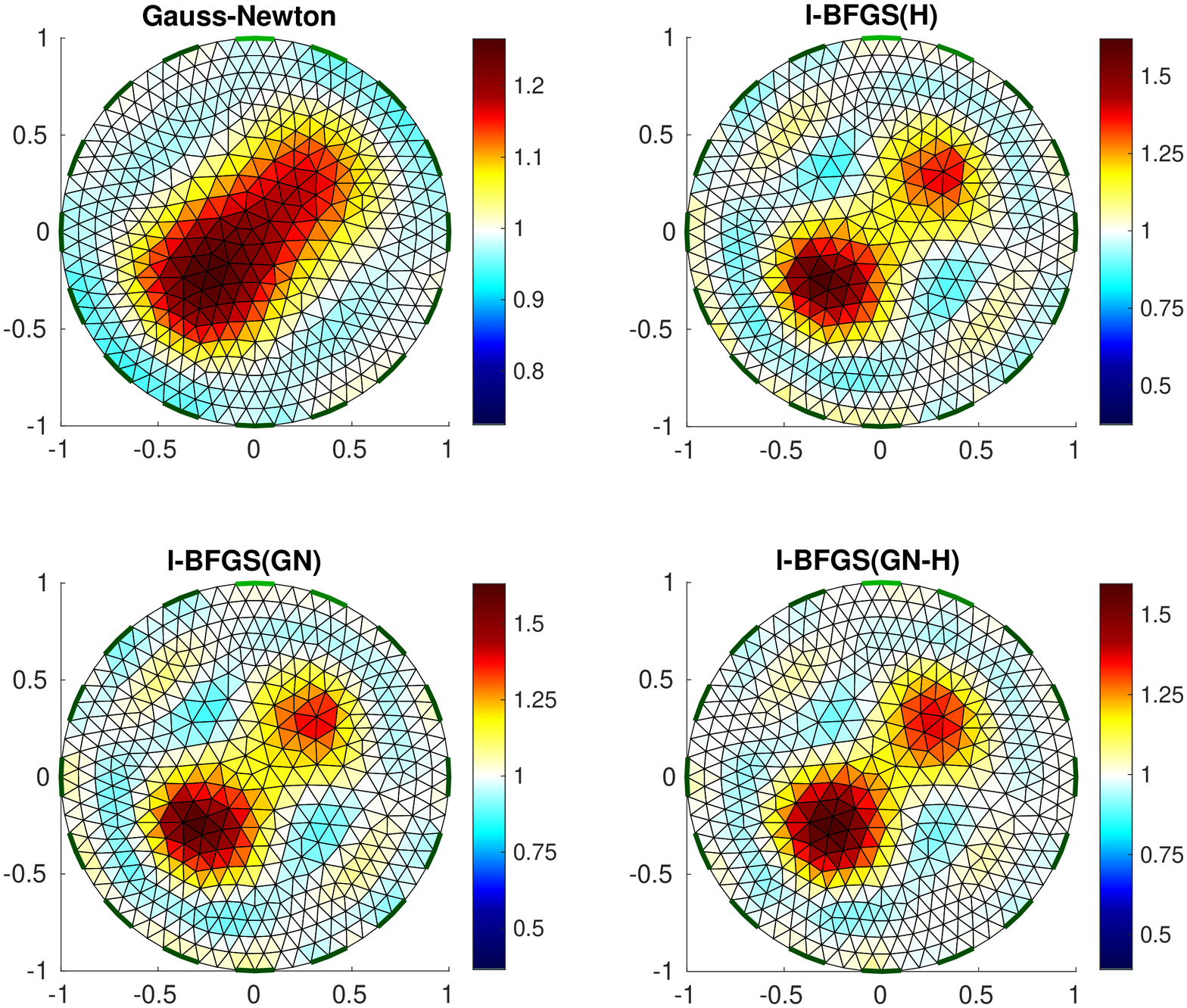}}}
	\caption{Final reconstruction results after the numerical stagnation condition (\ref{eq: numerical stagnation}) was met for a unit circular domain containing two inclusions of radius 0.16\(r_0\), separation 0.25\(r_0\).}
	\label{fig: recon 11}	
\end{figure}

\begin{figure}
	\centering
	\subfloat[True domain (left), and slice through \(x=y\) of reconstructions (right).]{\resizebox*{0.95\textwidth}{!}{
			\includegraphics[width= 0.75\textwidth]{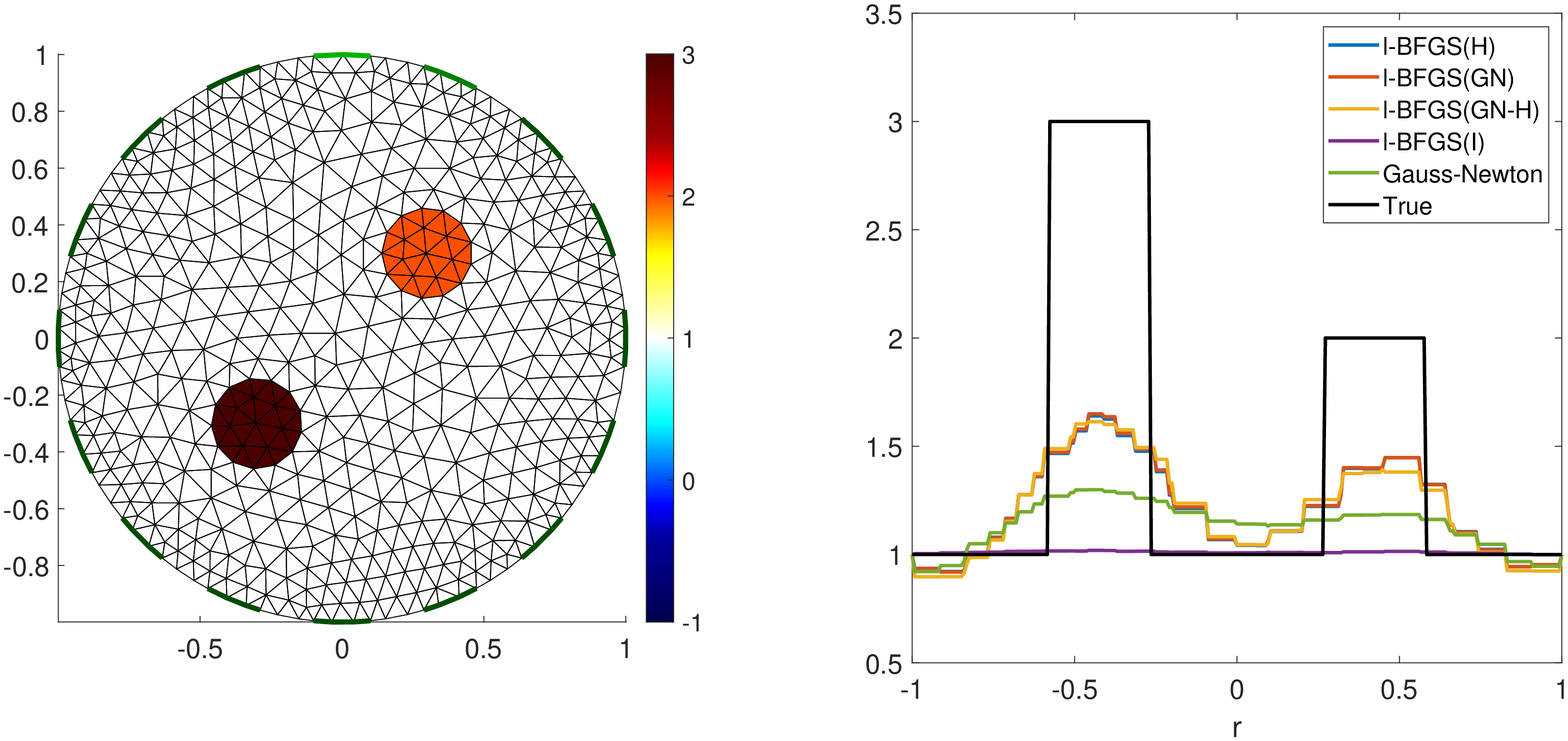}}}\\
	\subfloat[Reconstruction results with Gauss-Newton and l-BFGS]{\resizebox*{0.95\textwidth}{!}{	
			\includegraphics[width= 0.95\textwidth]{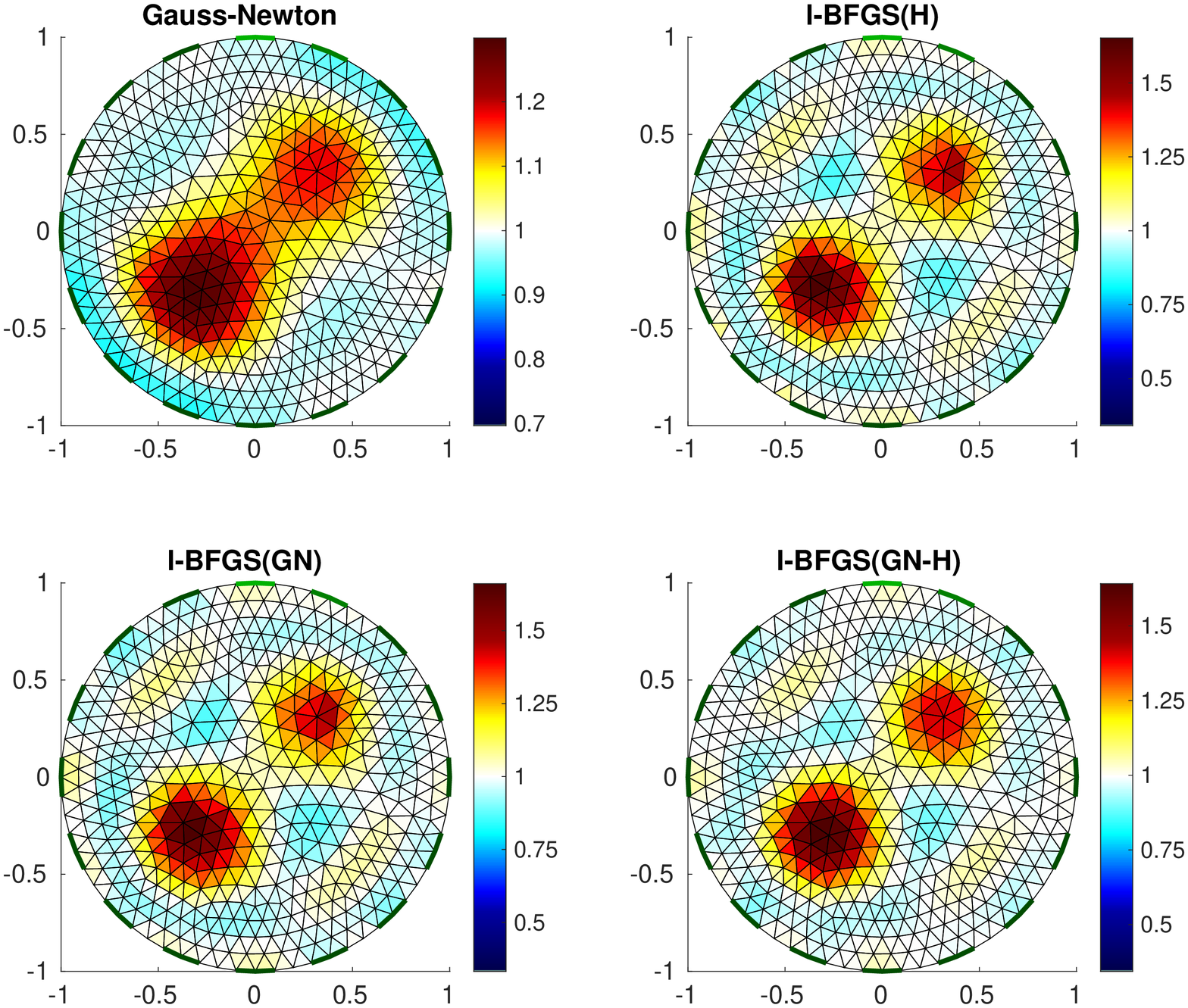}}}
	\caption{Final reconstruction results after the numerical stagnation condition (\ref{eq: numerical stagnation}) was met for a unit circular domain containing two inclusions of radius 0.16\(r_0\), separation 0.3\(r_0\).}
	\label{fig: recon 21}	
\end{figure}

\begin{figure}
	\centering
	\subfloat[True domain (left), and slice through \(x=y\) of reconstructions (right).]{\resizebox*{0.95\textwidth}{!}{
			\includegraphics[width= 0.75\textwidth]{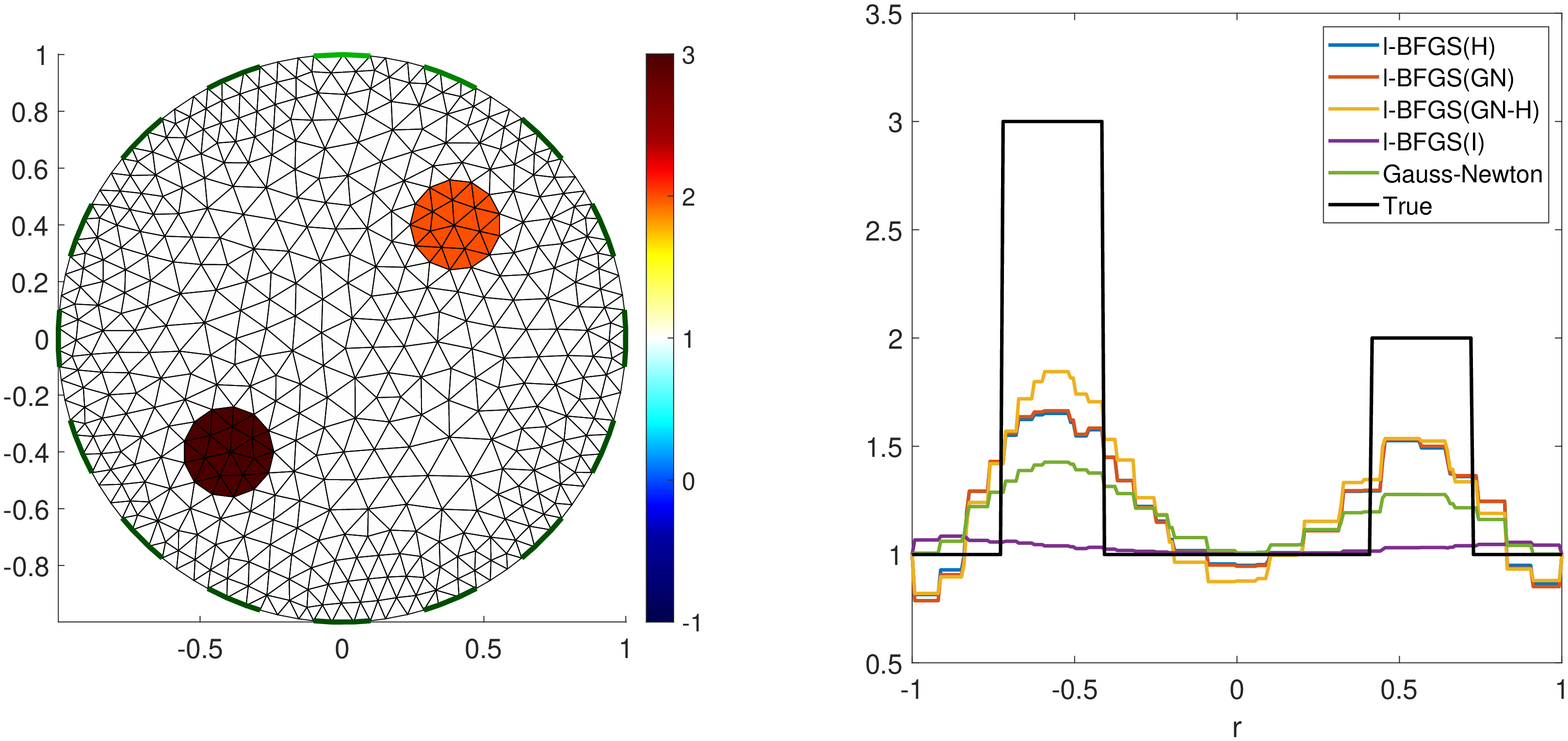}}}\\
	\subfloat[Reconstruction results with Gauss-Newton and l-BFGS]{\resizebox*{0.95\textwidth}{!}{	
			\includegraphics[width= 0.95\textwidth]{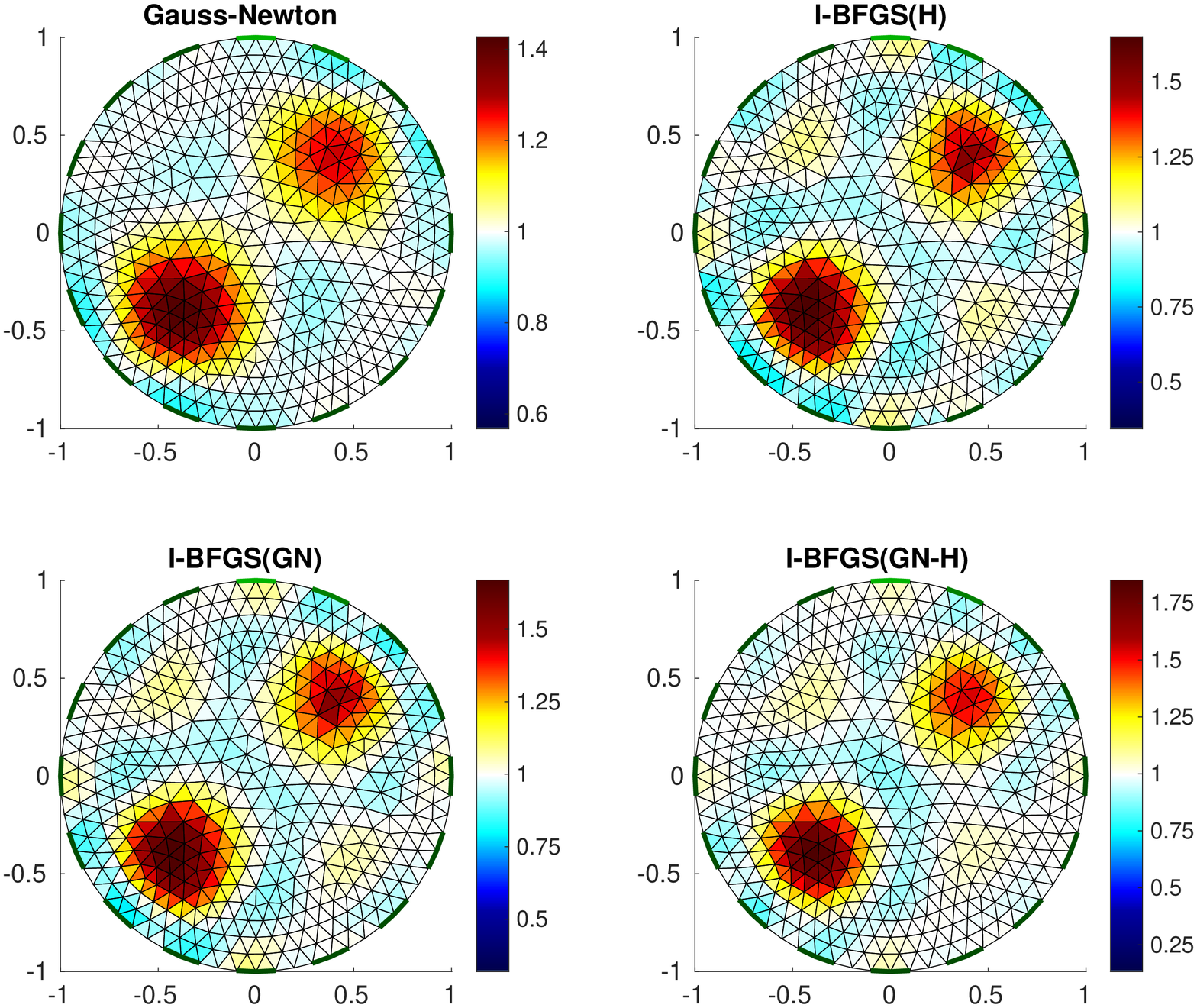}}}
	\caption{Final reconstruction results after the numerical stagnation condition (\ref{eq: numerical stagnation}) was met for a unit circular domain containing two inclusions of radius 0.16\(r_0\), separation 0.4\(r_0\).}
	\label{fig: recon 31}	
\end{figure}

\begin{figure}
	\centering
	\subfloat[True domain (left), and slice through \(x=y\) of reconstructions (right).]{\resizebox*{0.95\textwidth}{!}{
			\includegraphics[width= 0.75\textwidth]{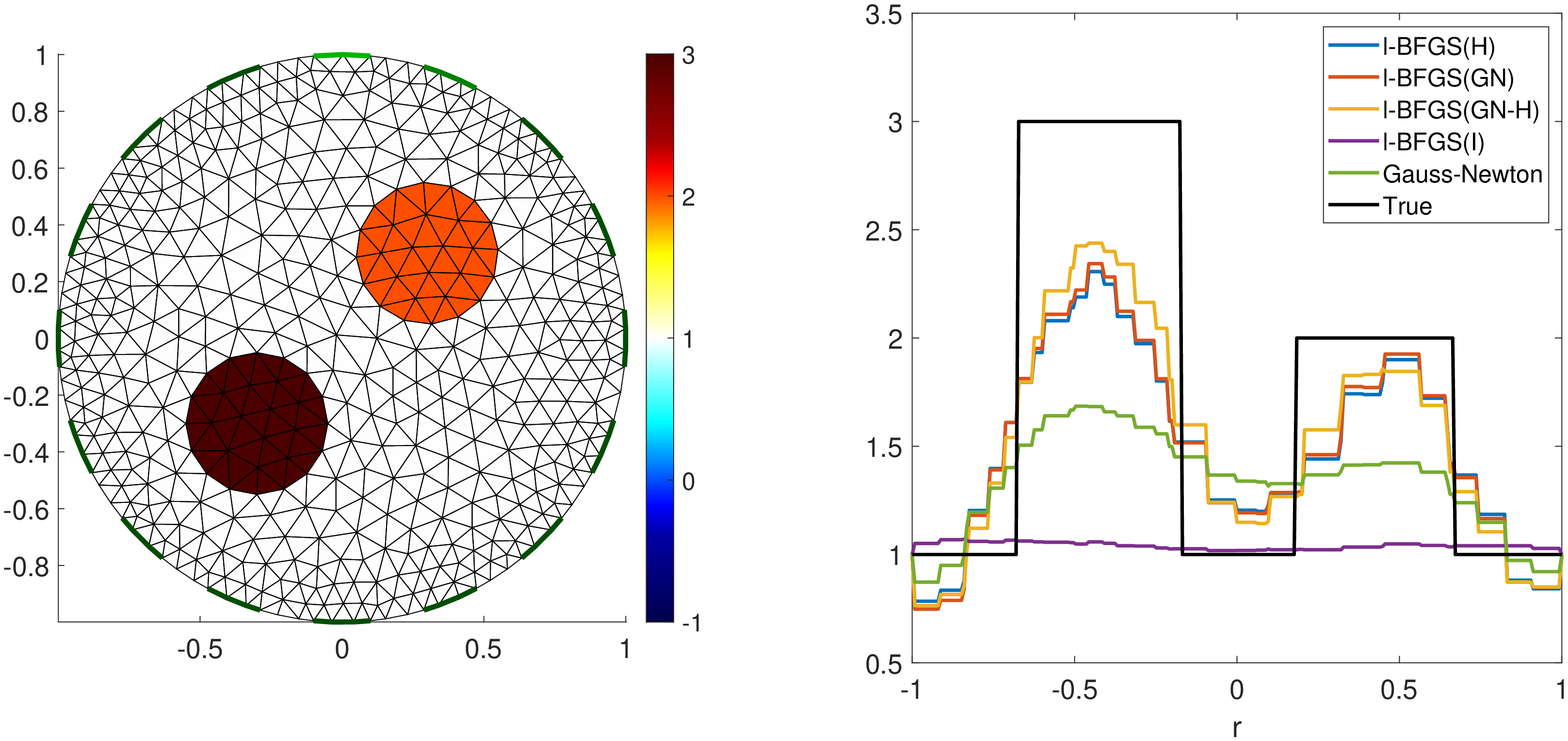}}}\\
	\subfloat[Reconstruction results with Gauss-Newton and l-BFGS]{\resizebox*{0.95\textwidth}{!}{	
			\includegraphics[width= 0.95\textwidth]{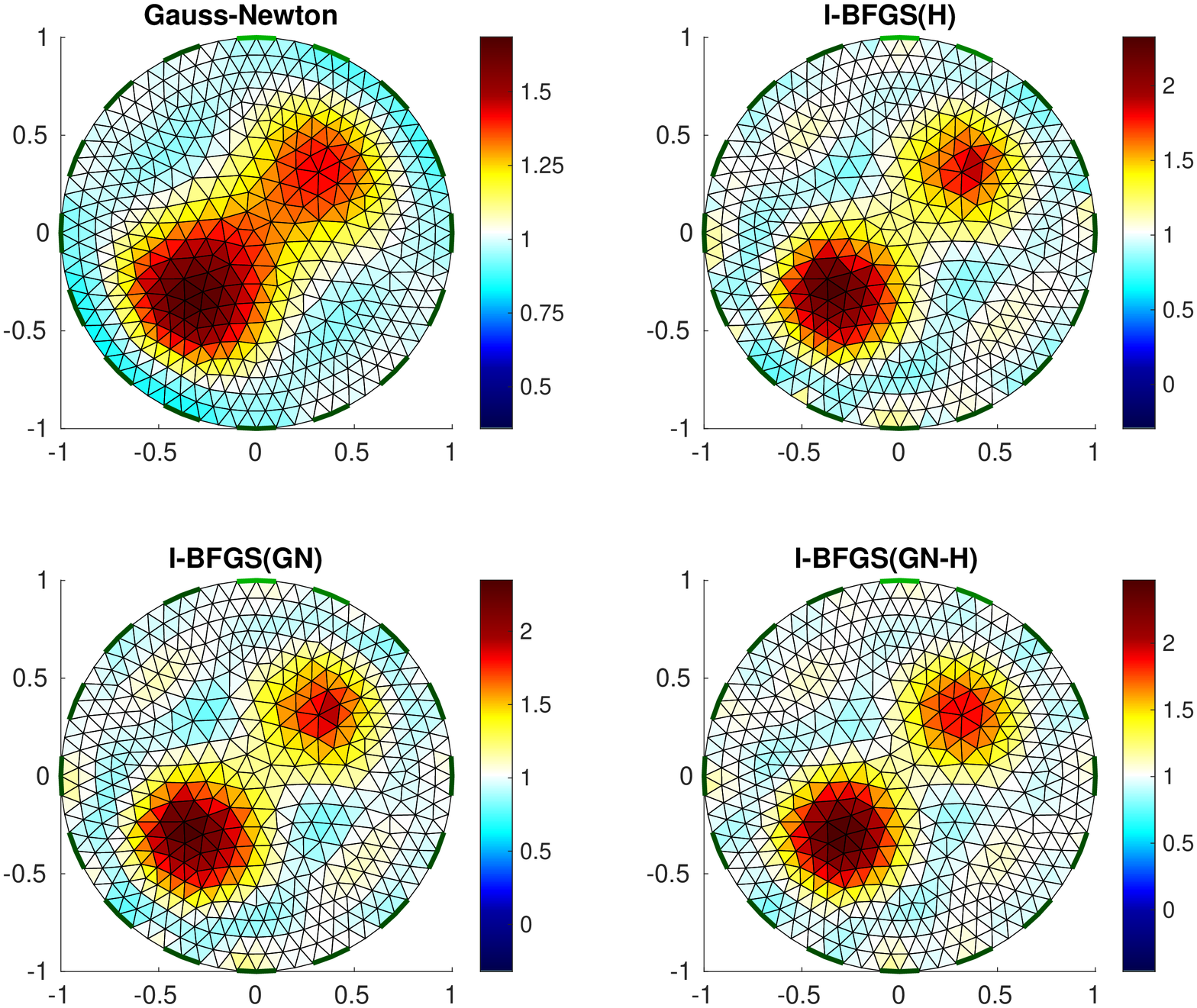}}}
	\caption{Final reconstruction results after the numerical stagnation condition (\ref{eq: numerical stagnation}) was met for a unit circular domain containing two inclusions of radius 0.25\(r_0\), separation 0.3\(r_0\).}
	\label{fig: recon 22}	
\end{figure}

The approximate diagonal Hessian is now used as the initial Hessian in an l-BFGS reconstruction scheme.  Four data sets are simulated for unit circle domains containing two circular inclusions with \(\sigma=2\sigma_0\) and \(\sigma=3\sigma_0\), respectively. For each, independent and identically distributed Gaussian noise is added with a signal-to-noise ratio of 50. Figures~\ref{fig: recon 11},~\ref{fig: recon 21} and \ref{fig: recon 31} show the true domain and reconstruction results for inclusions of radius \(0.16r_0\), separated from the origin by \(0.25r_0, 0.3r_0\) and \(0.4r_0\), respectively.  \figurename~\ref{fig: recon 22} shows the true domain and reconstruction results for inclusions of radius \(0.25r_0\), separated from the origin by \(0.3r_0\).  In each case, the reconstruction domain has \(N_m = N_E = 800\) elements.

Homogeneous backgrounds are used for easier assessment of reconstruction results versus objects size and distance apart, rather than the objects varying contrast to the background. However, our previous work for GPR imaging demonstrates the method can equally be applied to reconstruction problems with an unknown inhomogeneous background\cite{Watson15thesis}.

Four reconstructions are presented for each simulated data set, resulting from different reconstruction algorithms.  These are:
\begin{itemize}
	\item Gauss-Newton;
	\item {l-BFGS} with initial Hessian approximation given by
	\begin{itemize} 
		\item Polarization tensor approximate Hessian \(\tilde{\lmat{H}}\), referred to as l-BFGS(H), 
		\item \(\diag{(\lmat{J}^T\lmat{J})}\), referred to as l-BFGS(GN),
		\item \(\diag{(\lmat{J}^T\lmat{J})}\) plus the second derivative part of the polarization tensor approximate Hessian, referred to as l-BFGS(GN-H).
	\end{itemize}
\end{itemize}
The regularisation term for each was
\begin{equation}
R(\lvec{m}) = \|\lmat{L}\lvec{m}\|_2^2,
\end{equation}
with \(\lmat{L}\) the discrete Laplace operator in 2D, and a regularisation parameter of \(\lambda=5\times10^{-5}\) (chosen heuristically to provide the best stable Gauss-Newton results as a benchmark).  The regularisation term was also included in the initial Hessian approximation for each of the l-BFGS methods.  Each reconstruction was stopped when the relative change in residual met the numerical stagnation condition
\begin{equation}
\frac{\mathcal{J}(\lvec{m}^{[k-1]}) - \mathcal{J}(\lvec{m}^{[k]})}{\mathcal{J}(\lvec{m}^{[0]})} < 10^{-4},
\label{eq: numerical stagnation}
\end{equation}
which is the default condition of the EIDORS Gauss-Newton solver.

From figures~\ref{fig: recon 11} through \ref{fig: recon 31}, we see that the  variants of l-BFGS are generally better able to separate the two inclusions than Gauss-Newton, which tends to result in reconstructions with the two inclusions blurred together more.  The results of l-BFGS(H) and l-BFGS(GN) are almost indistinguishable, providing acceptable separation and sizing of the inclusions as well as showing some contrast between the two.  This is a favourable result, as l-BFGS(H) may be slightly less computationally expensive if the Jacobian matrix is not being calculated and stored (e.g. if an adjoint method is being used to calculate gradients).

In Figures~\ref{fig: recon 31} and \ref{fig: recon 22}, l-BFGS(GN-H) also results in a slightly better contrast estimation of the two inclusions, but otherwise similar (or indistinguishable) results.  Since this is also a computationally inexpensive addition to calculating the Gauss-Newton diagonal, we would also suggest this is a favourable result.  

\begin{figure}
	\centering
	\subfloat[Inclusion radius \(0.16r_0\), separation 0.25\(r_0\)]{\resizebox*{0.4\textwidth}{!}{
			\includegraphics[width=0.95\textwidth]{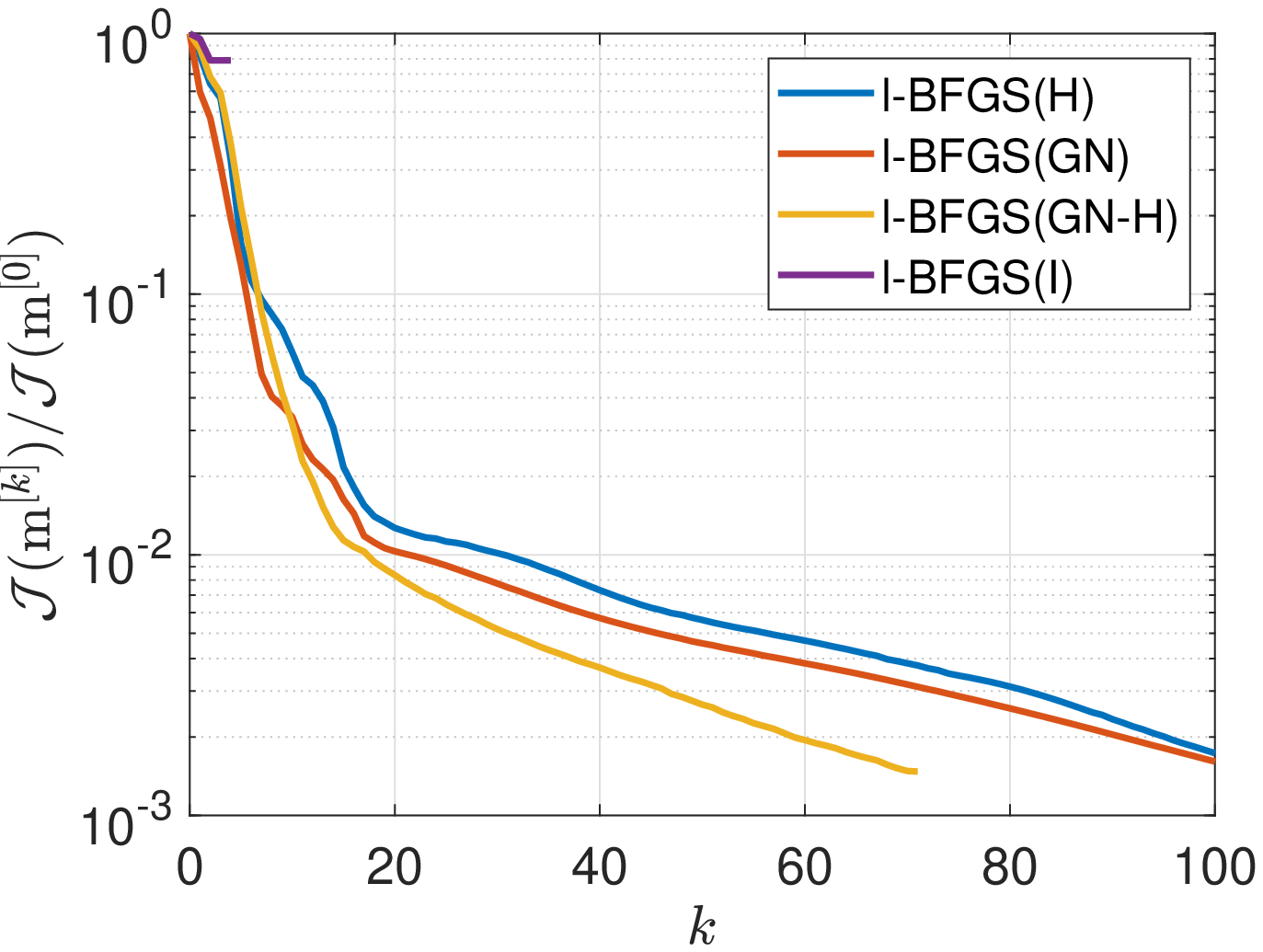}}}~
	\subfloat[Inclusion radius 0.16\(r_0\), separation 0.3\(r_0\)]{\resizebox*{0.4\textwidth}{!}{
			\includegraphics[width=0.95\textwidth]{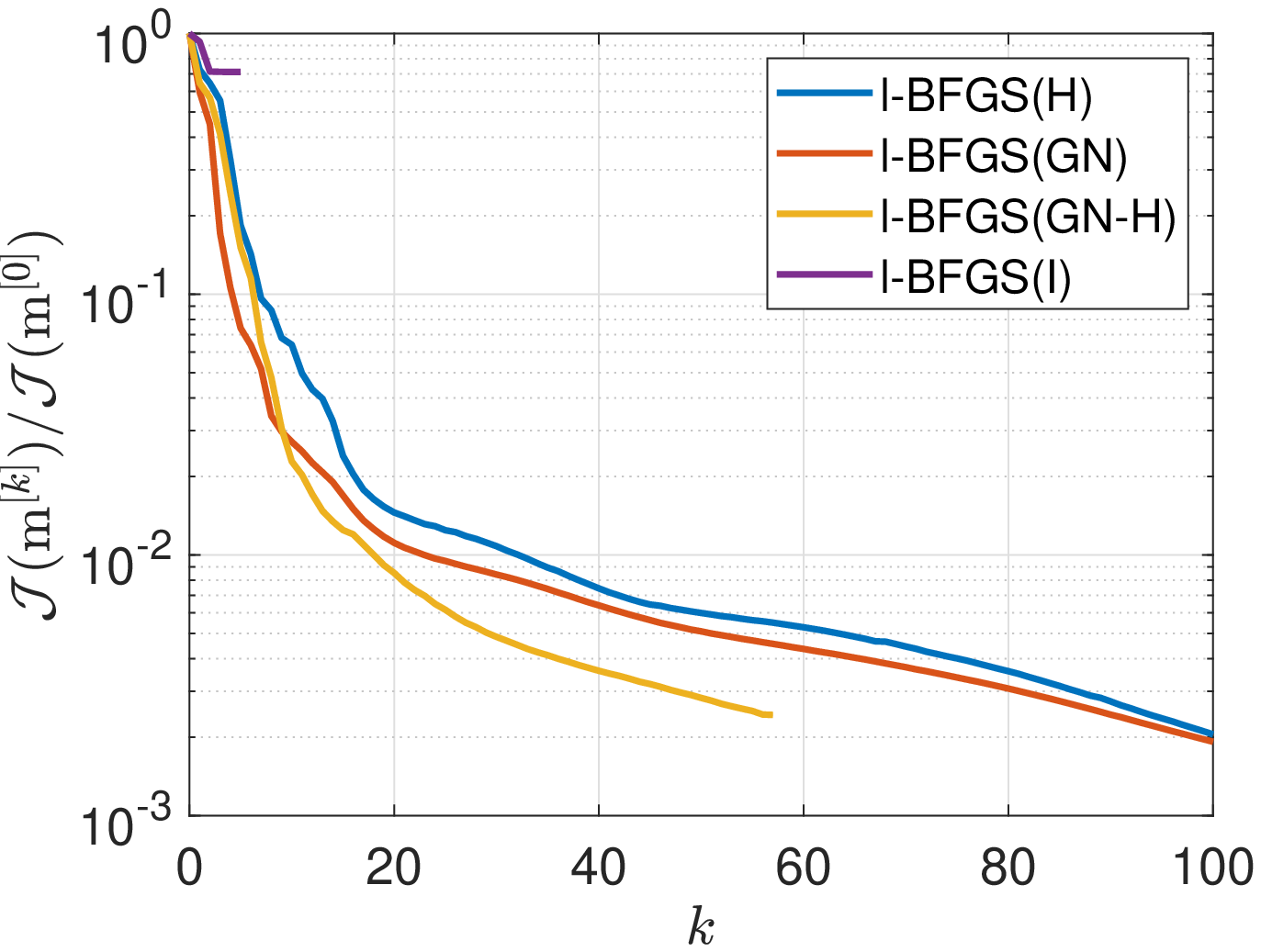}}}\\
	\subfloat[Inclusion radius 0.16\(r_0\), separation 0.4\(r_0\)]{\resizebox*{0.4\textwidth}{!}{
			\includegraphics[width=0.95\textwidth]{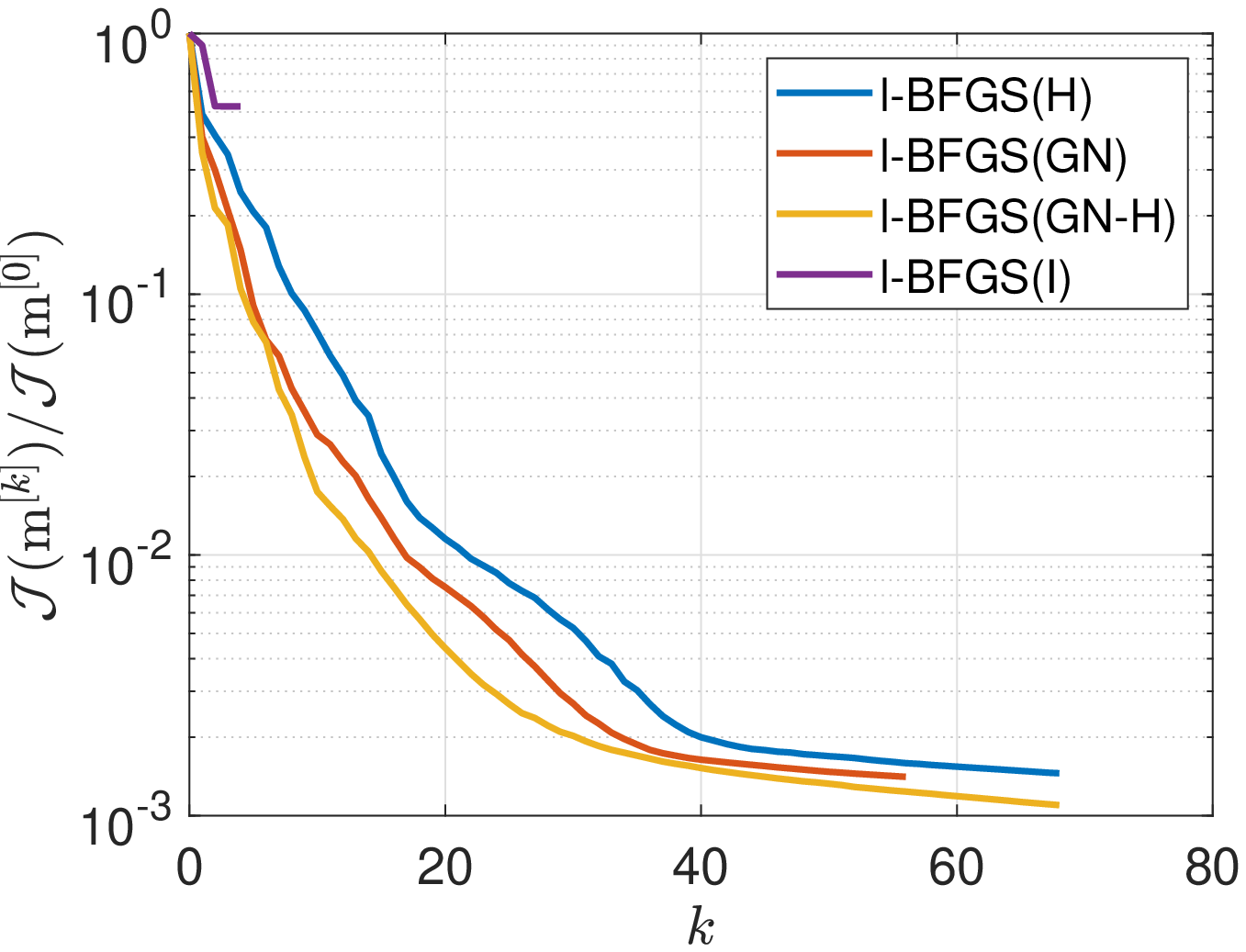}}}~
	\subfloat[Inclusion radius 0.25\(r_0\), separation 0.3\(r_0\)]{\resizebox*{0.4\textwidth}{!}{
			\includegraphics[width=0.95\textwidth]{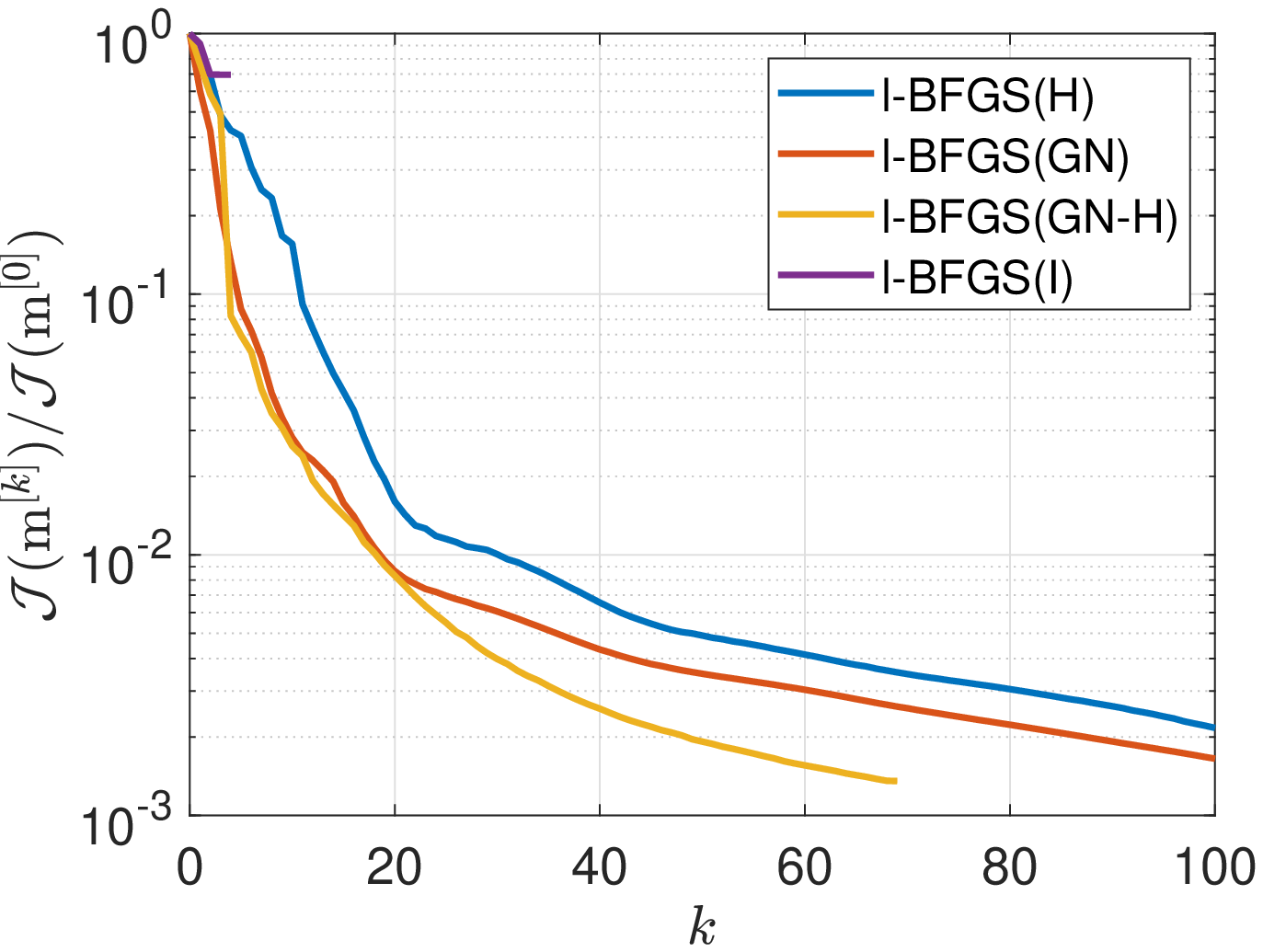}}}
	\caption{Objective function relative residual \(\mathcal{J}(\lvec{m}^{[k]})/\mathcal{J}(\lvec{m}^{[0]})\) against l-BFGS iteration \(k\) for the four reconstruction problems.}
	\label{fig: convergence}
\end{figure}

Since the visual quality of reconstruction results alone is not a sufficient measure of utility, the relative residual of the objective function against l-BFGS iteration is given in \figurename~\ref{fig: convergence}.  We see that l-BFGS(H) has a slightly slower rate of convergence than l-BFGS(GN), though the difference is generally marginal, and l-BFGS(GN-H) generally outperforms both. Indeed, for the three reconstruction problems with inclusions further from the boundary, l-BFGS(GN-H) converges in approximately 30 iterations fewer than the other two methods.   These results are displayed against l-BFGS iteration, and it may be interesting also to directly compare e.g. CPU time.  However, since EIDORS uses efficient factorisation and caching methods to form the Jacobian during the solution to the forward problem it would be impractical to develop our proposed method to a similar standard for a useful comparison.  Moreover, considering the additional cost of forming \(\tilde{\lmat{H}}\) is similar to that of \(\lmat{J}^T\lmat{J}\) (see Section~\ref{sec: quasi-newton}), we would expect little difference to the appearance of these results.

These results are not unexpected: l-BFGS(GN-H) provides the small second derivative correction to the exact Gauss-Newton diagonal, resulting in a more accurate approximation to the true Hessian and with that faster convergence.  With inclusions close to the boundary though, the polarization tensor approximation is known to break down, and the performance improvement is negligible (though it still continues to progress a few iterations beyond the other methods for a slightly improved solution).  For l-BFGS(H), it would appear that this second derivative correction has less impact than the inaccuracy of the first derivative approximation, hence performance is not improved over l-BFGS(GN).  

We also attempted to use l-BFGS initiated with the identity (denoted l-BFGS(I) in \figurename~\ref{fig: convergence}), but in each case it stagnated before any reasonable progress was made.  This is seen in the slices through reconstructions  in figures~\ref{fig: recon 11} through \ref{fig: recon 31} (purple line) as well as the convergence results in \figurename~\ref{fig: convergence}.  So, at the very least, l-BFGS(H) provides a computationally viable quasi-Newton method where memory constraints are a concern.  Then, l-BFGS(GN-H) provides a particularly computationally cheap and effective way to incorporate second derivative information as compared to calculating these terms directly.

\section{Conclusion}
Our motivation in this work has been to effectively deal with non-linearity in inverse boundary value problems for PDEs.  This includes both a lack of superposition of data caused by multiple nearby inclusions, as well as the saturation effect as the material contrast of a single inclusion increases.  In many problems involving well-isolated inclusions, solutions to the forward problem
exist in the form of asymptotic expansions with polarization tensors.  These can be used to find efficient and robust solutions to the inverse problem.  However, they become invalid where there are nearby inclusions and the lack of superposition becomes a more dominant effect.  In such cases, one often poses the inverse problem as a least-squares one, for which the Hessian matrix is known to be important in finding efficient and robust solutions.  

Our contribution has been to use the polarization tensors to derive a computationally cheap approximate diagonal Hessian matrix, which describes the saturation effect and ill-posedness due to parameter illumination.  We have proposed to use this either as an initial Hessian approximation in quasi-Newton schemes, or to provide the second derivative correction in addition to \(\diag(\lmat{J}^T\lmat{J})\) for such methods. Through our numerical experimentation with the EIT problem, we have shown that in the first case this provides a computationally viable quasi-Newton method for non-linear inverse problems which may have memory constraints. In the latter case, where one is able to calculate the Jacobian matrix itself, this provides a particularly inexpensive way to incorporate second derivative information which in some cases may have a significant positive impact on performance. Our numerical experimentation consisted of not just the reconstruction problems with simulated data, but also looked at the accuracy of derivative approximations themselves.  We present this as a proof-of-principle for this mixed-model approach to solving non-linear inverse problems, potentially usable for any imaging modality in which a polarization tensor approximation to the forward model exists. 

We have suggested some further work throughout this paper.  Firstly, our choice of ``closest fitting ellipse" to the triangular elements (in the electrical sense of (\ref{eq: explicit ellipse})) as the Steiner inellipse is taken for computational simplicity.  The development of either explicit results or a more precise heuristic rule might be used to improve this method.  Secondly, as only a proof-of-principle has been presented, the methods should be tested on realistic large 3D problems in which memory is a constraint, including other imaging modalities besides EIT.

\section*{Acknowledgement}
We would like to thank EPSRC for their support with grants EP/R002177/1, EP/L019108/1, EP/K00428X/1; the Royal Society for a Wolfson Research Merit Award; the Sir Bobby Charlton Foundation (formerly Find A Better Way); and the Defence, Science and Technology Laboratory.

\FloatBarrier
\bibliographystyle{tfnlm}
\bibliography{pbib}

\end{document}